\documentclass[10pt,a4paper]{amsart} 

\usepackage[english]{babel}
\usepackage[utf8]{inputenc}  

\usepackage{epsfig}
\usepackage{float}
\usepackage{subfigure}
\usepackage{moreverb}

\usepackage{xcolor}

\usepackage{afterpage}		
\usepackage{latexsym}
\usepackage{amsmath}
\usepackage{amssymb}

\usepackage[numbers,sort&compress]{natbib}

\renewcommand{\appendix}{
\setcounter{section}{0}
\renewcommand{\thesection}{\Roman{section}}
\vspace{0.5cm}
{\Large{\bf APPENDIX}}}

\newenvironment{remark}[1][Remark.]{\begin{trivlist}
\item[\hskip \labelsep {\bfseries #1}]}{\end{trivlist}}

\newcommand{\uu}{{\bf u}}

\newcommand{\ww}{{\bf w}}

\begin{document}

\title[Stable numerical discretisation of a new Boussinesq model]{A novel energy-bounded Boussinesq model and a well balanced and stable numerical discretisation}

\author{ Magnus Sv\"ard and Henrik Kalisch }

\date{\today}

\begin{abstract}
In this work, a novel Boussinesq system is put forward. The system is naturally
nonlinearly entropy/energy-stable, and is designed for problems with sharply 
varying bathymetric features. The system is flexible and allows tuning of the dispersive parameters 
to the relevant wavenumber range of the problem at hand. 
We present a few such parameter sets, including one that tracks the dispersive relation 
of the underlying Euler equations up to a nondimensional wavenumber of about $30$.

In the one-dimensional case, we design a stable finite-volume scheme and demonstrate 
its robustness and accuracy in a suite of test problems including Dingemans's wave experiment.
We generalise the system to the two-dimensional case and sketch how the numerical scheme 
can be straightforwardly generalised. 
\end{abstract}

\maketitle

\section{Introduction}


Boussinesq systems are generally coupled Partial Differential Equations (PDEs) that 
model the propagation of waves at the surface of a fluid under the assumption that 
the wavelength is long and the amplitude is small when compared to the depth. 
While a full description of surface waves should be given in terms of the Euler 
or Navier Stokes equations with free-surface boundary conditions \cite{Lannes2013},
the computational cost of integrating the full Euler equations on oceanic or even coastal scales 
is still prohibitive. In consequence, 
Boussinesq systems are highly relevant as simplified models 
for use in practical situations in coastal modeling.
The first system of this type appeared in the work of Boussinesq 
and was given on a flat bed \cite{Boussinesq1872}.
In \cite{Peregrine67}, the system was generalised to allow variations in the fluid bed,
and from that point on, a large number of of structurally similar Boussinesq-type 
systems have been derived.  

As mathematical models, the original Boussinesq system as well as many of its variants 
can be viewed as perturbations of the well-known shallow water equations, 
augmented with third-derivative terms that improve approximation with regard 
to the dispersive properties of surface waves as described by the full Euler equations.
One such system was given in \cite{MadsenSorensen92}. This system was derived under the assumption 
of slowly varying smooth bathymetry, and the resulting system is linearly stable. 
In other words, smooth perturbations of smooth solutions (water at rest) remain bounded.

If the bathymetry deviates too much from the above assumption, one cannot expect solutions 
to be bounded. Indeed, many Boussinesq systems have been shown to suffer from strong instabilities 
in cases where the bathymetry faetures sharp variations (see \cite{LovholtPedersen09} for example). 
To address this problem, considerable efforts have been made first and foremost by introducing
various forms of numerical dissipation with the hope of suppressing the instabilities.
However, as was also pointed out in \cite{LovholtPedersen09}, 
the instabilities are not only caused by the numerical discretisations, but may be inherent in 
the mathematical structure of the equations.
While some improvements can be made on the linear stability of the equations
\cite{simarro2015linear},
many Boussinesq systems are ill-posed outside of the linear regime.  
Needless to say, there is no ``numerical fix'' for instabilities observed in ill-posed equations.

Some Boussinesq systems are known to be well posed.
For example, a variant of the original Boussinesq system with a flat bed
was proven to be well posed in \cite{Schonbeck81}. Some other examples
pertaining to systems with slowly varying bathymetry and higher-order nonlinearities
are given in \cite{Lannes2013}, but in general the systems that accurately model the
dispersion relation such as the Nwogu system \cite{Nwogu1993} 
and the two-layer systems \cite{beji1996formal} are not known to be well posed.

Recently, in \cite{MadsenFuhrman20} so-called trough instabilities were discussed. 
As the name suggests, these instabilities appear in wave troughs
where high-frequency disturbances appear in certain situations.
The authors of \cite{MadsenFuhrman20} showed that models with linear dispersive terms 
(if linearly stable) do not suffer from these instabilities, 
while those with nonlinear dispersive terms often do.

At this point, it should be pointed out that well-posedness is usually not an integral part 
of the procedure used to derive new Boussinesq models. In many cases, these models are obtained 
by truncating a series expansion of the velocity potential in a similar fashion as used
in the work of Peregrine \cite{Peregrine67}, and subsequently adding other terms based 
on physical arguments. Well-posedness would thus have to verified once the model is defined. 

Proving well-posedness is generally very difficult for nonlinear partial differential equations.
However, a minimal requirement for well-posedness is that the solution remains appropriately bounded,
and in many cases, physical reasoning provides guidance when deriving a priori bounds of the solution. 
To make progress, we turn back to the full Euler equations for guidance. The Euler equations are based
on conservation of mass and momemtum, and these properties essentially carry over to the Boussinesq
system considered here (however the same is not true for all Boussinesq systems \cite{AliKalisch12}).
From the Euler system, the energy balance can be derived and one can conclude that energy 
may be converted from kinetic to potential energy and vice versa, and energy may also be dissipated,
i.e., turned into heat. Dissipated energy may not return to the system.

In accordance with the above considerations, a well posed Boussinesq-type system should at least
have bounded mechanical energy (in mathematics, such a condition is known as an \emph{entropy condition}). 
Contrary to this principle, we will demonstrate that common Boussinesq systems do \emph{not} 
feature bounded mechanical energy and thereby verify that they are ill-posed. This is bound to be an issue when non-linear effects are important. In consequence, it is futile to try to stabilise 
numerical approximation schemes; there will always be cases where they are 
unable to produce accurate results.

The main idea of this paper is to integrate the fundamental physical property 
of energy boundedness into the derivation process of a novel Boussinesq-type system. 
Indeed, we will show how that there is a set of dispersive terms which can be added to 
the shallow-water system such that the resulting system retains the property of energy conservation for a flat bathymetry.  
Once we have established a general system that admits bounded solutions, 
we determine the exact coefficients of the dispersive terms by requiring that 
solutions of the linearised problem have dispersive properties 
that are close to the full water-wave problem.

Next, we extend the new system to a varying bathymetry in an energy stable way. 
We make sure that the approximate dispersion relation is obtained at any constant depth. 
The latter property also dispenses with the necessity to switch from a Boussinesq model 
to the shallow-water model as the depth decreases. 
Such a device is used in some cases \cite{roeber2010shock,tonelli2011simulation},
but depends on accurate pinpointing of wave breaking \cite{bjorkavaag2011wave,boettger2020energetics}.
In our formulation this transition is built into the system as the dispersive terms 
essentially vanish if the depth becomes shallow enough. 
While some recent work has suggested how to handle wave breaking more efficiently
in Boussinesq-type systems \cite{Kazolea2018}, it is not
the purpose of the present work to investigate energy dissipation due to wave breaking.

The resulting system, being energy bounded, is subsequently discretised 
using the well-established framework of entropy-stable schemes, 
see \cite{Tadmor03, Carpenter_etal16, Fjordholm_etal11}, extended to the new system.

Finally, we validate the model against Dingemans's 
wave experiment (\cite{Dingemans94,Dingemans97}) and demonstrate its robustness 
for bathymetries with sharp features.

\section{The shallow water system}

We begin by introducing some notation. Let the $x$-axis be the horizontal dimension and $z$-axis the vertical. $P$ denotes the depth-integrated momentum, and $d$ denotes the local depth, i.e, from the bottom to the actual surface (not the still water depth). Furthermore, $v=P/d$ is the depth integrated velocity; $b(x)$ is the bathymetry, i.e., a function that defines the bottom in relation to some reference value on the $z$-axis. $H$ is the height to the still-water surface from the same reference point. Hence, $d+b-H$ is the deviation from the still-water surface (surface elevation).
As usual, $g$ is the local gravitational constant. 

In these variables, the 1-D shallow-water equations are given by:
\begin{align}
  d_t + P_x &=0, \label{SWE} \quad x\in \Omega,\quad t>0, \\
  P_t + \left(\frac{P^2}{d}\right)_x + gd(d+b-H)_x&=0. \nonumber
\end{align}
where $\Omega$ denotes the spatial domain. The entropy-pair of (\ref{SWE}) is given by:
\begin{align}
  U(\uu)&=\frac{1}{2}\frac{P^2}{d}+\frac{1}{2}gd^2+gdb,\label{U} \\
  F(\uu)&=\frac{1}{2}\frac{P^3}{d^2}+gdP+gPb,\label{F}
\end{align}
where $\uu=(d,P)^T$. For the shallow-water equations, the entropy $U(\uu)$ is also the total mechanical energy.  As discussed above, in the absence of external energy sources and if there is no energy flux through the boundaries, solutions of (\ref{SWE}) must satisfy the additional constraint that the total energy in the system does not exceed its initial value. That is, kinetic and potential energy may dissipate into heat, but heat may not transform back into mechanical energy.
\begin{remark}
Although $U(\uu)$ is the mechanical energy and has no relation to the specific entropy in thermodynamics, we shall call it the \emph{entropy function}, or \emph{entropy} for short as is common in the mathematical and numerical literature on conservation laws. However, when relating our results to physics we will sometimes use the term \emph{mechanical energy}. Finally, we use the term \emph{energy} for the $L^2$-norm of perturbations in the linear analysis. 
\end{remark}

From this observation, the entropy balance can be derived. To this end, we introduce the gradient of the entropy with respect to the conservative (principle) variables, which is termed  \emph{entropy variables}. Here, they are:
\begin{align}
\ww^T=(g(d+b)-\frac{P^2}{2d^2},\,\frac{P}{d}).\label{ent_var}
\end{align}
The entropy balance is obtained by contracting (\ref{SWE})  with (\ref{ent_var}).
\begin{align}
  0&=\left(g(d+b)-\frac{P^2}{2d^2}\right)\left(d_t+P_x\right)
  + \frac{P}{d}\left(P_t + \left(\frac{P^2}{d}\right)_x + gd(d+b)_x\right).\label{ent1}
\end{align}
We need the identities:
\begin{align}
  -\frac{P^2}{2d^2}d_t+\frac{P}{d}P_t&=\partial_t\left(\frac{P^2}{2d}\right)\nonumber \\
  g(d+b)d_t&=\partial_t\left(\frac{1}{2}gd^2+gdb\right)\nonumber \\
  \left(g(d+b)-\frac{P^2}{2d^2}\right)P_x  + 
\frac{P}{d}\left(\left(\frac{P^2}{d}\right)_x + gd(d+b)_x\right)&=\partial_x\left(gP(d+b)+\frac{1}{2}\frac{P^3}{d^2} \right). \label{aux1}
\end{align}
By introducing $v=P/d$, the momentum part follows from the following calculation,
\begin{align}
  -\frac{1}{2}v^2(dv)_x+v(dv^2)_x=
  -\frac{1}{2}v^2(dv)_x+\frac{1}{2}v(dv^2)_x+\frac{1}{2}(v^2(dv)_x+v^2dv_x)=\nonumber \\
    \frac{1}{2}v(dv^2)_x+\frac{1}{2}v^2dv_x=\frac{1}{2}(dv^3)_x=\partial_x\left(\frac{P^3}{2d^2}\right).\label{aux2}
\end{align}
Inserting (\ref{aux1}) and (\ref{aux2}) in (\ref{ent1}) results in,
\begin{align}
  0&=\partial_t U(\uu) +\partial_x F(\uu),\label{ent2}
\end{align}
which is the entropy equality. Next, we assume that the system does not interact with its surroundings, i.e., $v=0$ at the left and right boundaries. Upon integration in time and space, (\ref{ent2}) gives the entropy bound,
\begin{align}
\int_{\Omega}U(\uu(\cdot,t))\,dx \leq \int_{\Omega}U(\uu(\cdot,0))\,dx.
\end{align}
That is, the mechanical energy is bounded by its intial state. More generally, with appropriate boundary conditions, the mechanical energy cannot grow unboundedly.

\begin{remark}
Note that the entropy $U$ is allowed to grow due to entropy that is entering through the boundaries. For a problem where there is no interaction with the surroundings, the system should not allow $U$ to grow.
  \end{remark}

\section{The system of Madsen and S\o rensen}

Next, we consider the entropy balance of Boussinesq-type systems. To this end, we use the system proposed 
in \cite{MadsenSorensen92}, but we emphasise that it is merely used as a prototypical example. 
The system in \cite{MadsenSorensen92} takes the form.
\begin{align}
  d_t + P_x &=0,\nonumber \\
  P_t + \left(\frac{P^2}{d}\right)_x + gd(d+b-H)_x+\Phi_B&=0,\label{eq1}
\end{align}
where
\begin{align}
\Phi_B=-(B+1/3)P_{xxt}-Bgh_s^3(d+b-H)_{xxx},\label{BMS}
\end{align}
are the dispersive Boussinesq terms. Furthermore, $h_s$ is the still water depth and $B$ a constant.

The entropy analysis differs from (\ref{SWE}) only in the  $\Phi_B$ terms. The result is,
\begin{align}
  0&=\partial_t U(\uu) +\partial_x F(\uu)
  + \frac{P}{d}\left(-(B+1/3)P_{xxt}-Bgh_s^3(d+b)_{xxx}\right).\label{ent2}
\end{align}
For the last two terms to be entropy consistent, they need to form complete derivatives and/or positive quadratic terms, which they do not. Hence, they remain indefinite and may (depending on the particular flow state) cause an unbounded growth of $U(\uu)$.

We remark that the system (\ref{eq1}) is linearly stable.
Hence, it works for small amplitude waves and smoothly and slowly varying bathymetry.
Indeed those are the assumptions in the derivation of (\ref{eq1}).
However, exactly how small the amplitude and bathymetry variation have to be,
is impossible to quantify a priori. 
At any rate, in order to bound solutions beyond the linear regime,
a non-linear bound is indispensible.
\emph{In fact, the lack of a non-linear entropy bound for (\ref{eq1}) (and similar models)
  explains the lack of stability when the bathymetry is rough.}

\subsection{Other Boussinesq approximations}

 The starting point when deriving (\ref{eq1}) was the system derived in (\cite{Peregrine67}) where
\begin{align}
\Phi_B=\frac{1}{6}h^3(\frac{P}{h})_{xxt}-\frac{1}{2}h^2(P)_{xxt}.\label{peregrine}
\end{align}
Interestingly, these dispersive terms do not lead to a bounded entropy either. Contracting with the entropy variables gives the following contribution to the entropy balance.
\begin{align*}
\frac{P}{d}\Phi_B=\frac{P}{d}\left(\frac{1}{6}h^3(\frac{P}{h})_{xxt}-\frac{1}{2}h^2(P)_{xxt}\right).
\end{align*}
As these terms do not form complete derivatives and positive quadratic terms, the system (\ref{eq1}) with (\ref{peregrine}) does not bound the entropy either. Hence, it is not surprising that systems derived from (\ref{peregrine}) may be unstable.

\section{Entropy-bounded Boussinesq system}

Having pointed to the mathematical issues with current Boussinesq systems, we will follow a new procedure for deriving a new model:
\begin{enumerate}
\item The system has to be a dispersive perturbation of (\ref{SWE}). That is, it inherits the entropy of
  the shallow-water system.
\item The dispersive terms are chosen to ensure entropy boundedness and each term is scaled with a coefficient.
\item The entropy analysis provides constraints on the coefficients. 
\item We choose the coefficients by matching the dispersive relation of the Boussinesq system with that of the Euler equations for a flat bathymetry. We propose a few different options leading to systems of varying complexity and accuracy.
\item We recast the system to allow varying bathymetry while making sure that
  \begin{itemize}
    \item entropy boundedness is uncompromised,
  \item  it leads to the optimised dispersion relation at any constant depth, and
  \item the system reduces to (\ref{SWE}) as the depth goes to zero.
  \end{itemize}
\end{enumerate}

\begin{remark}
  The last property will dispense with the necessity to switch from a Boussinesq approximation
  to the shallow-water system near shores.
\end{remark}

We observe that Boussinesq systems generally feature dispersive terms in the momentum and the surface elevation, i.e., $S=(d+b-H)$. However, as we noted above, momentum dispersion is not compatible with entropy boundedness. Hence, we consider
\begin{align}
  d_t+P_x&= (\alpha(d+b)_{xx})_x \label{eqSane1}\\
  P_t + \left(\frac{P^2}{d}\right)_x + gd(d+b-H)_x &=(\alpha v(d+b)_{xx})_x+ \beta \left(\frac{P}{d}\right)_{xxt}+\gamma\left(\frac{P}{d}\right)_{xxx},\label{eqSane2}
\end{align}
where $\alpha,\beta,\gamma$ are constants.

As before, we demonstrate entropy boundedness by contracting (\ref{eqSane1})-(\ref{eqSane2}) with the entropy variables (\ref{ent_var}). The shallow-water part is the same as before and we obtain,
\begin{align}
  \partial_t U(\uu) +\partial_x F(\uu)=T_\alpha+T_\beta+T_\gamma,\label{ent_balance}
\end{align}
where $T_{\alpha,\beta,\gamma}$ are the contributions from the dispersive terms. Clearly, they must not induce an unbounded growth and we calculate their contributions  one-by-one:
\begin{align*}
  T_\alpha = (g(d+b)-\frac{P^2}{2d^2})(\alpha(d+b)_{xx})_x+(\alpha v(d+b)_{xx})_x(\frac{P}{d}).
\end{align*}
We recast the first term as,
\begin{align*}
  g(d+b)(\alpha(d+b)_{xx})_x=
  g((d+b)(\alpha(d+b)_{xx}))_x-g(d+b)_x\alpha(d+b)_{xx}&=\\
  g((d+b)(\alpha(d+b)_{xx}))_x-\frac{1}{2} \alpha g(\frac{(d+b)_x^2}{2})_x,
\end{align*}
and the second
\begin{align*}
  -\frac{v^2}{2}(\alpha(d+b)_{xx})_x+v(\alpha v(d+b)_{xx})_x&=\\
  (-\frac{v^2}{2}\alpha(d+b)_{xx})_x  +vv_x\alpha(d+b)_{xx})
  +(v\alpha v(d+b)_{xx})_x  -v_x\alpha v(d+b)_{xx}&=\\
  (\frac{v^2}{2}\alpha(d+b)_{xx})_x. 
\end{align*}
where we have used $P/d=v$. Then,
\begin{align*}
 T_\alpha= g((d+b)(\alpha(d+b)_{xx}))_x-\frac{1}{2} \alpha g(\frac{(d+b)_x^2}{2})_x+ (\frac{v^2}{2}\alpha(d+b)_{xx})_x. 
\end{align*}
Clearly, $T_\alpha$ does not violate entropy boundedness as it is expressed in divergence form. Furthermore, the constant $\alpha$ can be chosen freely without violating entropy boundedness.

Next, we consider the system
\begin{align*}
T_{\gamma}=v\gamma v_{xxx}= (\gamma vv_{xx})_x-\frac{\gamma}{2}(v_x^2)_x.
\end{align*}
As above, $T_\gamma$ has been recast to divergence form and does not contribute to entropy growth.

Turning to $T_\beta$:
\begin{align*}
T_\beta=v\beta v_{xxt} = (\beta v v_{xt})_x- \beta v_xv_{xt}= (\beta v v_{xt})_x- \frac{\beta}{2} (v_x^2)_t.
\end{align*}
Here, we obtain one term in divergence form and a temporal derivative of $v_x^2$ 
that requires further attention.  
Upon integration of (\ref{ent_balance}) in time and space (domain $(0,1)$), 
we have the entropy balance
\begin{align*}
\int_0^T \int_0^1(U(\uu)_t + F(\uu)_x)\,dx\,dt = \int_{0,0}^{T,1} (T_\alpha+T_\beta+T_\gamma)\,dx\,dt.
\end{align*}
For simplicity, we consider the periodic case in which case all boundary terms cancel 
and the above relation turns into
\begin{align*}
\int_0^T \int_0^1 (U(\uu)_t\,dx\,dt=  \int_{0,0}^{T,1} (T_\beta)\,dx\,dt= - \int_{0,0}^{T,1} ( \frac{\beta}{2} (v_x^2)_t)\,dx\,dt, \nonumber
\end{align*}
or,
\begin{align}
\int_{0}^{1} \left(U(\uu)|_{t=T}-U(\uu)|_{t=0}+\frac{\beta}{2} (v_x^2)|_{t=T}-\frac{\beta}{2} (v_x^2)|_{t=0}\right)\,dx&=0.\label{ent_balance}
\end{align}
Naturally, we assume that initial data is sufficiently bounded: $\|U(\uu)\|_1|_{t=0}\leq Constant$ and $\|v_x\|_2|_{t=0}\leq Constant$. Then, if we assume that $\beta>0$, we obtain a bound on the positive quantity $\int_{0}^{1}(U(\uu)|_{t=T}+\frac{\beta}{2} (v_x^2)|_{t=T})\,dx\leq Constant$. 

\begin{remark}
Note that the structure of the dispersive terms in (\ref{eqSane1})-(\ref{eqSane2}) is a requirement for stability. 
In particular, it is not possible to add a dispersive term for the surface elevation 
in one equation without a balancing term in the other. 
\end{remark}
We summarise our findings so far: 
The Boussinesq system (\ref{eqSane1}) and (\ref{eqSane2}) has a bounded entropy, i.e., mechanical energy, when
\begin{itemize}
\item $\alpha,\gamma$ are arbitrary constants, and
\item  $\beta$ is an arbitrary \emph{non-negative} constant.
\end{itemize}
That is, we have satisfied the three first items in our list.

\subsection{Linear dispersion relation}

To obtain concrete models, we will use the dispersive relation for linear waves in order to choose the coefficients $\alpha,\beta,\gamma$. (Item 4 in the list.)

To this end, we linearise  (\ref{eqSane1})-(\ref{eqSane2}) around $v=0$ and $d=d_0$, and assume that $b=constant$ . Let $d=d_0+s$, where $s$ is the surface elevation, i.e., the perturbation of the depth variable $d$. Furthermore, we denote the velocity perturbation as $v'$. Linearizing (\ref{eqSane1}), results in
\begin{align}
  (d_0+s)_t+((0+v')(d_0+s))_x&= (\alpha((d_0+s+b_0)_{xx})_x,\quad \textrm{or}\nonumber \\
  s_t+d_0v'_x&= \alpha s_{xxx}   \label{eqLin1}
\end{align}
where a quadratically small term has been omitted in (\ref{eqLin1}). In the same way, the linearisation of (\ref{eqSane2}) yields,
\begin{align}
  (v'd_0)_t + gd_0s_x &= \beta \left(v'\right)_{xxt}+\gamma\left(v'\right)_{xxx}.\label{eqLin2}
\end{align}
Next, we recast (\ref{eqLin1}) and (\ref{eqLin2}) as a single second-order wave equation by differentiating (\ref{eqLin1}) by $t$ and (\ref{eqLin2}) by $x$
\begin{align}
  s_{tt}+d_0v'_{xt}&= \alpha s_{xxxt}, \label{stt}\\
  d_0v'_{xt} + gd_0s_{xx} &= \beta \left(v'\right)_{xxxt}+\gamma\left(v'\right)_{xxxx},\nonumber
\end{align}
 and combine the two,
\begin{align*}
  d_0 (\alpha s_{xxxt} -s_{tt})d_0^{-1} + gd_0s_{xx} &= \beta \left( (\alpha s_{xxxt} -s_{tt})d_0^{-1}\right)_{xx}+\gamma\left(v'\right)_{xxxx}.
\end{align*}
To obtain an equation in $s$ alone, we differentiate in time
\begin{align*}
  \alpha s_{xxxtt} -s_{ttt} + gd_0s_{xxt} &= \beta \left( (\alpha s_{xxxt} -s_{tt})d_0^{-1}\right)_{xxt}+\gamma\left(v'\right)_{xxxxt},
\end{align*}
and use (\ref{stt}) again,
\begin{align*}
  \alpha s_{xxxtt} -s_{ttt} + gd_0s_{xxt} &= \beta \left( (\alpha s_{xxxt} -s_{tt})d_0^{-1}\right)_{xxt}+\gamma d_0^{-1}\left(\alpha s_{xxxt}-s_{tt}\right)_{xxx}.
\end{align*}
All terms are differentiated with respect to time at least once, and we simplify to
\begin{align*}
  \alpha s_{xxxt} -s_{tt} + gd_0s_{xx} &= \beta \left( (\alpha s_{xxxt} -s_{tt})d_0^{-1}\right)_{xx}+\gamma d_0^{-1}\left(\alpha s_{xxx}-s_{t}\right)_{xxx},
\end{align*}
which is the linearised equation for constant bathymetry.
Next, we insert the wave solution $s=\exp(ikx-i\omega t)$.
\begin{align*}
  \alpha (-ik^3)(-i\omega) +\omega^2 + gd_0(-k^2) =&\, \beta \left( (\alpha (-ik^3)(-i\omega) +\omega^2)d_0^{-1}\right)(-k^2)\nonumber\\
 & +\gamma d_0^{-1}\left(\alpha (-ik)^3+i\omega\right)(-ik)^3.
\end{align*}
Simplify,
\begin{align*}
-\alpha k^3\omega+\omega^2-gd_0k^2=\beta\alpha d_0^{-1}(k^5\omega)-\beta d_0^{-1} k^2\omega^2-\gamma\alpha d_0^{-1} k^6 + \gamma d_0^{-1}\omega k^3.
\end{align*}
The characteristic equation is thus a quadratic equation in $\omega$:
\begin{align}
(1+\beta d_0^{-1}k^2)\omega^2 + (-\alpha - \beta \alpha d_0^{-1}k^2-\gamma d_0^{-1})k^3\omega -gd_0k^2+\gamma\alpha d_0^{-1} k^6=0\label{char}
\end{align}
One procedure to determine the coefficiencts $\alpha,\beta,\gamma$ is to solve (\ref{char}),
and make a polynomial approximation that is subsequently matched with a polynomial approximation
of the the dispersive target relation, obtained from the Euler equations \cite{Whitham74}, 
\begin{align}
\omega_{Euler}=\sqrt{gk}\sqrt{\tanh(d_0k)}\label{dispersion}.
\end{align}
We begin by making a polynomial approximation of (\ref{dispersion}). To this end, we need,
\begin{align*}
  \tanh(x)&=x-\frac{x^3}{3}+\frac{2x^5}{15}....\,\,\textrm{and,}\\
  \sqrt{1+x}&=1+\frac{x}{2}-\frac{x^2}{8}+\frac{x^3}{16}+...
\end{align*}
Then we recast (\ref{dispersion}) as
\begin{align}
  \omega \approx \pm\sqrt{gk}\sqrt{d_0k-\frac{(d_0k)^3}{3}+\frac{2}{15}(d_0k)^5\hdots}&=\nonumber\\
  \pm\sqrt{gd_0k^2}\sqrt{1-\frac{(d_0k)^2}{3}+\frac{2}{15}(d_0k)^4\hdots}&\approx\nonumber\\
  \pm k\sqrt{gd_0}(1-\frac{(d_0k)^2}{6}+\frac{1}{15}(d_0k)^4)&\label{approx_disp}.
\end{align}
Next, we solve (\ref{char}). We begin by only considering dispersion in the surface elevation $s$. That is, $\gamma=0,\beta=0$. Then (\ref{char}) simplifies to,
\begin{align}
\omega^2 -\alpha k^3\omega -gd_0k^2=0.\label{char_alpha}
\end{align}
Hence, by choosing
\begin{align}
  \alpha&=-\sqrt{gd_0}\frac{d_0^2}{3}\label{set1},\quad\textrm{(Set 1)} \\
  \beta&=\gamma=0,\nonumber
\end{align}
the two first terms of $\omega$ coincide with (\ref{approx_disp}).

\begin{figure}[h]
\centering
\subfigure{\includegraphics[width=8cm]{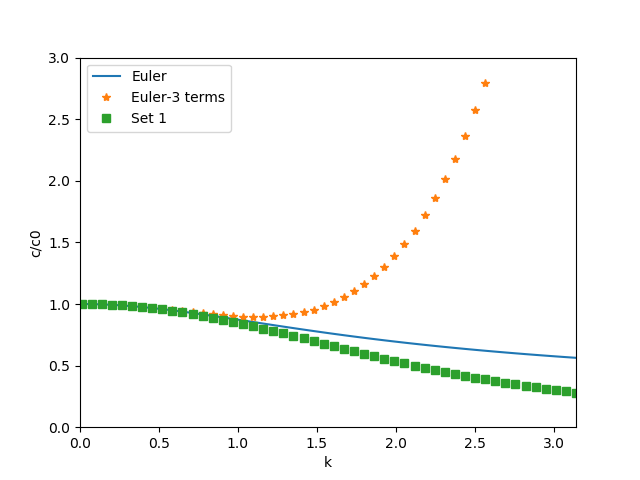}}
\subfigure{\includegraphics[width=8cm]{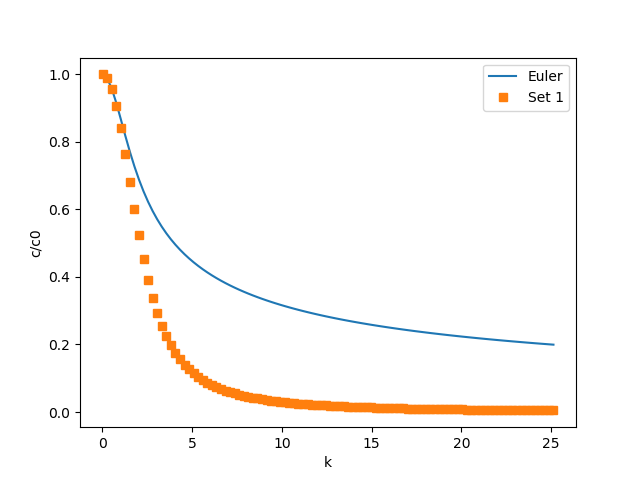}}
\subfigure{\includegraphics[width=8cm]{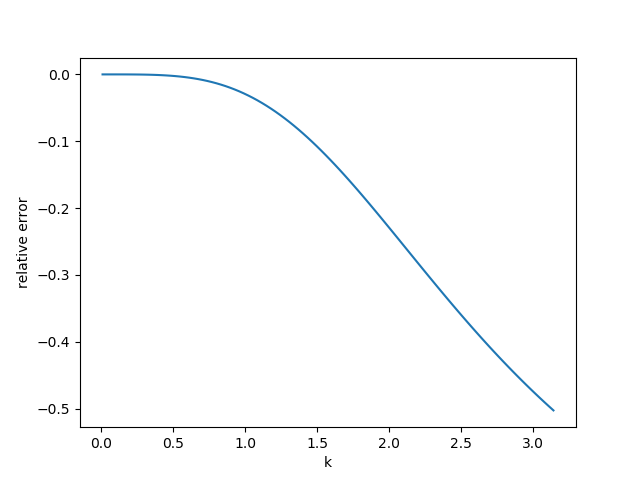}}
\caption{Dispersion relations and relative errors for ``Set 1'' (\ref{set1}).  ``Euler'' refers to (\ref{dispersion}) and ``Euler 3-terms'' to (\ref{approx_disp}).}\label{fig_disp_alpha}
\end{figure}

In Fig. \ref{fig_disp_alpha}, the dispersion relations are plotted (normalised with $c_0=\sqrt{g d_0}$). ``Euler'' refers to the dispersion relation (\ref{dispersion}) of the full Euler equations. ``Euler-3 terms'' displays the first three terms of the series expansion of the Euler relation given in (\ref{approx_disp}). ``Set 1'' refers to (\ref{set1}).  The relative error
\begin{align}
  error_r=\frac{|\omega_{approx}-\omega_{Euler}|}{\omega_{Euler}},\label{rel_err}
\end{align}
is fairly small up to $k=1$ (less than $4\%$). However, as seen in the Figure displaying (\ref{set1}) up to $k=8\pi$, the errors become large for higher frequencies. Hence, this simple model is only useful for long wavelengths. 

We also observe that (\ref{set1}) tracks (\ref{dispersion}) further than the approximate dispersion relation (\ref{approx_disp}). Hence, if the goal is to match (\ref{dispersion}) at high frequencies, one should \emph{not} proceed by matching more coefficients of (\ref{approx_disp}). Instead, we will attempt to choose coefficients that track (\ref{dispersion}) directly.
To this end, we make the system (\ref{eqSane1})-(\ref{eqSane2}) dimensionless by introducing
\begin{align}
\tilde d=\frac{d}{d_0},\,\tilde P=\frac{P}{c_0d_0},\,\tilde x=\frac{x}{d_0},\, \tilde t=\frac{t c_0}{d_0},
\end{align}
where $c_0=\sqrt{g d_0}$ and $d_0$ is a reference depth. (Since we assume a constant bathymetry, $d_0$ is the still water depth.) We obtain
\begin{align}
  \tilde d_{\tilde t}+\tilde P_{\tilde x}&= (\tilde\alpha(\tilde d+\tilde b)_{\tilde x\tilde x})_{\tilde x}, \label{nondimSane1}\\
 \tilde P_{\tilde t} + \left(\frac{\tilde P^2}{\tilde d}\right)_{\tilde x} + \tilde d(\tilde d+\tilde b)_x &=(\tilde \alpha \tilde v(\tilde d+\tilde b)_{\tilde x\tilde x})_{\tilde x}+ \tilde \beta \left(\frac{\tilde P}{\tilde d}\right)_{\tilde x\tilde  x\tilde t}+\tilde \gamma\left(\frac{\tilde P}{\tilde d}\right)_{\tilde x\tilde x\tilde x},\label{nondSane2}
\end{align}
where the dimensionless parameters, $\tilde \alpha,\tilde \beta,\tilde \gamma$ are
\begin{align}
  \alpha&=\tilde \alpha\sqrt{gd_0}d_0^2,\nonumber\\
  \beta&=\tilde \beta d_0^3,\nonumber\\
  \gamma&=\tilde \gamma\sqrt{gd_0}d_0^3.\nonumber
\end{align}
\begin{remark}
The effect of this non-dimensionalisation is that ``$g=1$'' and ``$d_0=1$'' in the characteristic equation (\ref{char}). 
\end{remark}
Next, we sweep the parameter space numerically in search of a set $\tilde \alpha, \tilde \beta, \tilde \gamma$ that gives a system with good dispersive properties.

\begin{figure}[h]
\centering
\subfigure{\includegraphics[width=6cm]{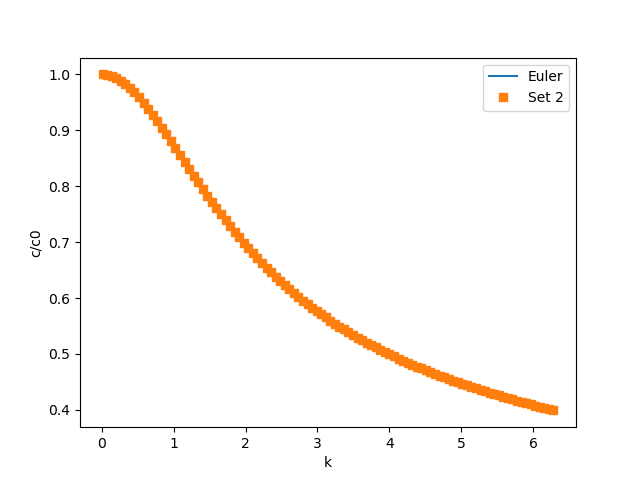}}
\subfigure{\includegraphics[width=6cm]{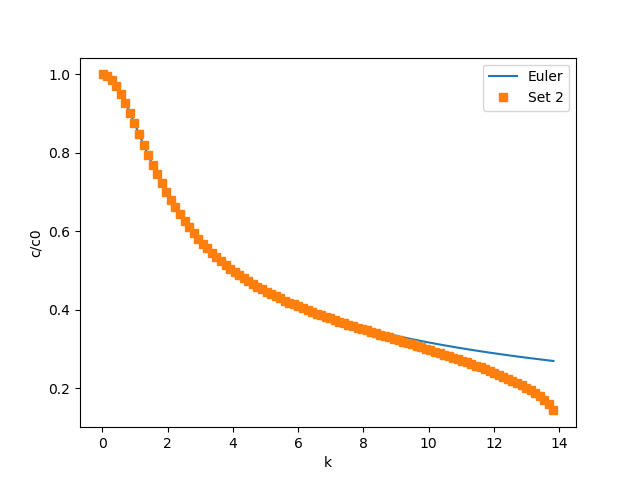}}
\caption{Dispersion relation for ``Set 2'' (\ref{set2}).  ``Euler'' refers to (\ref{dispersion}). The right panel displays (\ref{set2}) on $k=(0,4.4\pi)$.}\label{fig_set2}
\end{figure}

By trying to minimise the max-norm of the relative error (\ref{rel_err}) on the interval $k=(0,2\pi)$, we found,
\begin{align}
\tilde \alpha&= 0.0004040404040404049,\nonumber \\
\tilde \gamma&= 0.15707070707070708,\quad \textrm{(Set 2)}\label{set2}\\
\tilde \beta&= 0.49292929292929294.\nonumber 
\end{align}
The dispersion relation for (\ref{set2}) is depicted in Fig. \ref{fig_set2}. It lies on top of (\ref{dispersion}) and the maximal relative error is less than $0.9\%$.
Clearly, this is an excellent choice for $k$ up to $2\pi$ but the curve diverges for higher $k$.

\begin{figure}[h]
\centering
\subfigure{\includegraphics[width=6cm]{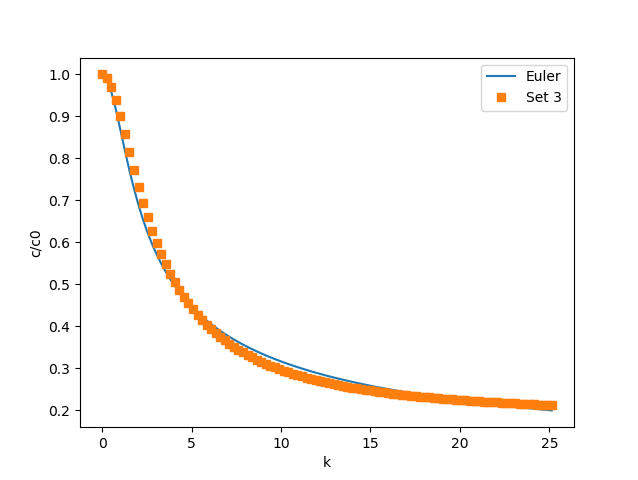}}
\subfigure{\includegraphics[width=6cm]{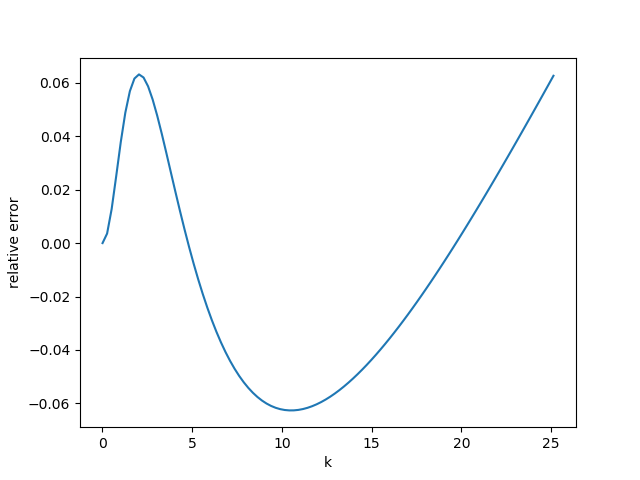}}
\caption{Dispersion relation for ``Set 3'' (\ref{set3}).  ``Euler'' refers to (\ref{dispersion}). The right panel displays the relative errors for (\ref{set3}).}\label{fig_set3}
\end{figure}

\begin{figure}[h]
\centering
\subfigure{\includegraphics[width=6cm]{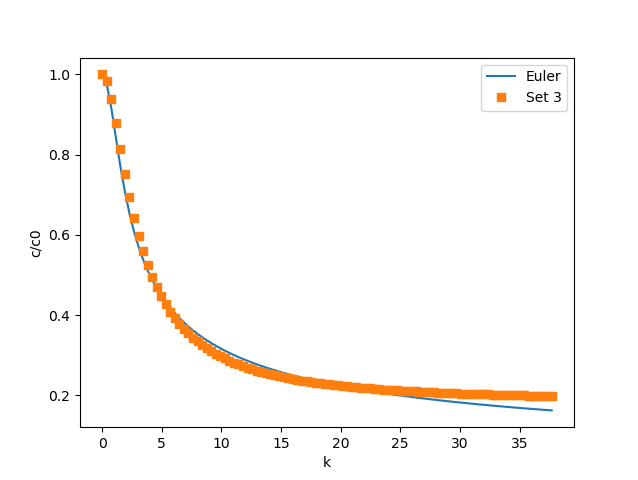}}
\caption{Dispersion relation for ``Set 3'' (\ref{set3}) on $k=(0,12\pi)$.  ``Euler'' refers to (\ref{dispersion}).}\label{fig_set3_high}
\end{figure}

In the next set, we have searched $k=(0,8\pi)$ for the smallest relative error. We found the following parameters:
\begin{align}
\tilde \alpha&=  0.0,\nonumber\\
\tilde \gamma&= 0.0521077694235589,\quad \textrm{(Set 3)}\label{set3}\\
\tilde \beta&= 0.27946992481203003.\nonumber 
\end{align}
The dispersion relation is shown in Fig. \ref{fig_set3}. Here, the match is visibly not perfect and the maximal relative error is less than $6.3\%$ on $k=(0,8\pi)$.
Furthermore, in Fig. \ref{fig_set3_high}, its properties are shown up to $k=12\pi$. The relative errors increase to $22\%$ as the dispersion relation deviates from (\ref{dispersion}).
As is the nature of polynomial approximations, a better match at the rapid decay of $c$ in the low-frequency range comes at the cost of a poorer approximation at high frequencies and vice versa. Hence, by sacrificing accuracy for $k\in (0,2\pi)$ one can match (\ref{dispersion}) better for higher frequencies than (\ref{set3}) does.

\begin{figure}[h]
\centering
\subfigure{\includegraphics[width=6cm]{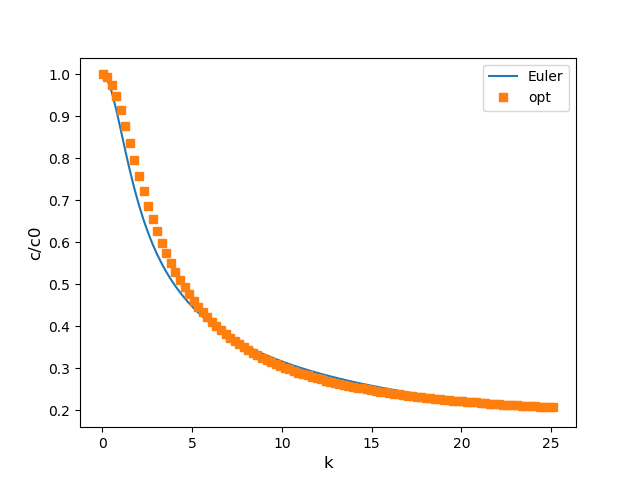}}
\subfigure{\includegraphics[width=6cm]{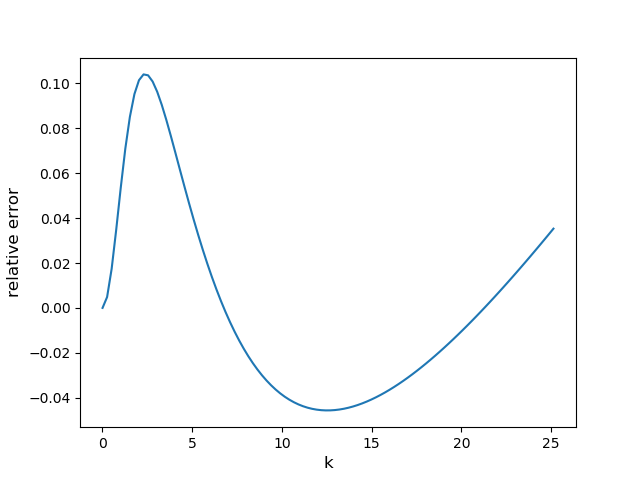}}
\caption{Dispersion relation for ``Set 4'' (\ref{set4}).  ``Euler'' refers to (\ref{dispersion}). The right panel displays the relative errors for (\ref{set4}).}\label{fig_set4}
\end{figure}

Moreoever, one can use weight functions when optimising the coefficients to improve the dispersive characteristics for certain frequencies. A simple such change is to use the absolute instead of the relative error when optimising the coefficients. That is, replace (\ref{rel_err}) with,
\begin{align}
  error_a=|\omega_{approx}-\omega_{Euler}|.\label{abs_err}
\end{align}
Looking for the minimal absolute error in $k=(0,2\pi)$ (with $\alpha=0$) results in,
\begin{align}
\tilde \alpha&=  0.0\nonumber\\
\tilde \gamma&= 0.04034343434343434,\quad \textrm{(Set 4)}\label{set4}\\
\tilde \beta&= 0.2308939393939394\nonumber
\end{align}
The relative error favours accuracy of low frequencies while (\ref{abs_err}) does not favour any range. The maximum relative error with this choice is  $10.4\%$. The dispersion relation is depicted in (\ref{fig_set4})
The examples (\ref{set1}),(\ref{set2}), (\ref{set3}) and  (\ref{set4}, are meant to demonstrate that, to a certain degree, it is possibly to adjust the system (\ref{eqSane1})-(\ref{eqSane2}) to ones needs. As long as $\beta\geq 0$, the system has a bounded mechanical energy.

\subsection{Varying bathymetry and shallow water}\label{sec:var_bat}

We have come to the 5th point in the list of requirements but before we address it, we make some general remarks.

When the water depth is small, it is common practice to ``turn off'' the Boussinesq terms and let the shallow-water
equations (\ref{SWE}) govern the flow \cite{roeber2010shock,bacigaluppi2020implementation}.
However, there are drawbacks with such an approach: 1) It is unclear what model, i.e. set of PDEs, that is being solved since that depends on some more or less arbitrary threshold in the code. Hence, it is nearly impossible to verify that the code solves any equation. 2) ``turning off terms'' in a numerical code is generally a highly unstable procedure. The reason is that these terms are now present on a subdomain which is not closed with boundary conditions, which is often an ill-posed procedure in itself. Sufficient artificial diffusion may control the instabilities on a given grid, but grid convergence is generally impossible to achieve. We emphasise that although grid convergence is often not carried out in engineering applications due to the lack of computational resources, it is the \emph{possibility} to grid refine to a converged solution, that is the sole guarantee that a numerical solution at hand is an approximation of the solution to the PDE at all.

Rather than having some threshold for shutting off the dispersive terms in a numerical solver, it is much more desirable to design the PDE such that the dispersive terms vanish at shallow water, as it enables both validation and verification, and not the least provides a possibility to reproduce results by other researchers.

To achieve this, we have to modify (\ref{eqSane1}) and (\ref{eqSane2}). We introduce a new set of model coefficients: $\hat \alpha(x),\hat \beta (x),\hat \gamma(x)$. (We will return to their specific forms.) Then we consider the following system:
\begin{align}
  d_t+P_x&= (\hat\alpha(\hat \alpha(d+b)_x)_{x})_x \label{eqShallow1}\\
  P_t + \left(\frac{P^2}{d}\right)_x + gd(d+b-H)_x &=(\hat \alpha v(\hat \alpha(d+b)_{x})_x)_x,\nonumber\\
  & + \left(\hat \beta \left( \frac{P}{d}\right)_x\right)_{xt}+\frac{1}{2}\left(\hat \gamma\left(\frac{P}{d}\right)_x\right)_{xx}
+\frac{1}{2}\left(\hat \gamma\left(\frac{P}{d}\right)_{xx}\right)_{x}.\nonumber \\
  \label{eqShallow2}
\end{align}

We denote the contributions to entropy balance, analogous to (\ref{ent_balance}), of the new modified dispersive terms, corresponding to $T_{\alpha,\beta,\gamma}$, as $M_{\alpha,\beta,\gamma}$. We have
\begin{align}
  M_\gamma=&\frac{1}{2}\left(v\left(\hat \gamma v_x\right)_{xx}
  +v\left(\hat \gamma v_{xx}\right)_{x}\right)=\nonumber\\
&\frac{1}{2}\left((v\left(\hat \gamma v_x\right)_{x})_x
+\left(v(\hat \gamma v_{xx})\right)_{x}-(v_x\hat\gamma v_x)_x
\right)\nonumber\\
  M_\beta=&v\left(\hat \beta\left( \frac{P}{d}\right)_x\right)_{xt}=(v\hat\beta ( v)_{x})_{xt}-\frac{1}{2}((\sqrt{\hat \beta} v_x)^2)_t,\nonumber
\end{align}
and finally,
\begin{align*}
  M_\alpha = (g(d+b)-\frac{P^2}{2d^2})(\hat\alpha(\hat\alpha(d+b)_{x})_x)_x+(\hat\alpha v(\hat\alpha(d+b)_{x})_x)_x(\frac{P}{d}).
\end{align*}
First, we handle
\begin{align*}
  g(d+b)(\hat\alpha(\hat\alpha(d+b)_{x})_x)_x&=\\
  g((d+b)\hat\alpha(\hat\alpha(d+b)_{x})_x)_x-  g((d+b)_x\hat\alpha)(\hat\alpha(d+b)_{x})_x)&=\\
  g((d+b)\hat\alpha(\hat\alpha(d+b)_{x})_x)_x
  -\frac{1}{2}  g(\frac{\hat\alpha(d+b)_x^2}{2})_x.
\end{align*}
As in the constant coefficient case, we obtain a divergence form. Next, we manipulate
\begin{align*}
  -\frac{v^2}{2}(\hat\alpha(\hat\alpha(d+b)_{x})_x)_x+v(\hat\alpha v(\hat\alpha(d+b)_{x})_x)_x&=\\
  (-\frac{v^2}{2}\hat\alpha(\hat\alpha(d+b)_{x})_x)_x+vv_x\hat\alpha(\hat\alpha(d+b)_{x})_x\\
  +(v\hat\alpha v(\hat\alpha(d+b)_{x})_x)_x-v_x\hat\alpha v(\hat\alpha(d+b)_{x})_x)
  &=\\
 (\frac{v^2}{2}\hat\alpha(\hat\alpha(d+b)_{x})_x)_x.
\end{align*}
Then,
\begin{align*}
 M_\alpha= g((d+b)\hat\alpha(\hat\alpha(d+b)_{x})_x)_x
  -\frac{1}{2}  g(\frac{\hat\alpha(d+b)_x^2}{2})_x + (\frac{v^2}{2}\hat\alpha(\hat\alpha(d+b)_{x})_x)_x
\end{align*}
The entropy estimate at time $T$ is obtained in the same way as the constant coefficient case. (C.f. (\ref{ent_balance}).)

Next, we turn to the choice of parameters and introduce the still water depth $h(x)$. The goal is to make the system reduce to (\ref{eqSane1})-(\ref{eqSane2}) for any constant depth. Clearly, that is ensured by the following choices
\begin{align}
  \hat \alpha^2&=\tilde \alpha\sqrt{gh(x)}h(x)^2,\nonumber\\
  \hat \beta&=\tilde \beta h(x)^3,\label{hat_par}\\
  \hat \gamma&=\tilde \gamma\sqrt{gh(x)}h(x)^3,\nonumber
\end{align}
where $\tilde\alpha,\tilde\beta,\tilde \gamma$ are the same coefficients as in the constant-depth system. For instance, (\ref{set1}),(\ref{set2}) and (\ref{set3}).

We also observe that, in addition to $\tilde \beta\geq 0$, we must now require that  $\tilde\alpha\geq 0$. (However, this does not seem to be a serious constraint since all attempts to optimise the three parameters have led to non-negative values of $\tilde \alpha$.)

Furthermore, if the depth $h(x)\rightarrow 0$, then the parameters (\ref{hat_par}) vanish. Thus, all three subitems in the 5th entry of our list of model constraints are satisfied.

\subsection{Well balanced and entropy-consistent artificial diffusion.}\label{sec:art_diff}
Although Boussinesq systems should not be diffusive, artificial diffusion terms are always required in order to stabilise numerical schemes in the presence of shocks (breaking waves) and to avoid inadmissible solutions (negative depths). Such terms also play a profound role when proving convergence of numerical approximations to non-linear conservation laws.

Diffusive terms are only stabilising if the they fit into the entropy framework of the problem. That is, if they diffuse the entropy $U(\uu)$. Furthermore, we require that artificial diffusion is well-balanced.

Hence, we will study what form of diffusion that may be added to the system (\ref{eqSane1})-(\ref{eqSane2}) (or equivalently (\ref{eqShallow1})-(\ref{eqShallow2})). To emphasise that this is a vanishing regularisation, we scale the diffusive terms with small, possibly vanishing, parameters $\epsilon>0,\delta>0$. Neglecting the dispersive terms, whose entropy consistent form we now know, we consider the following shallow-water system:
\begin{align}
  d_t+P_x&=(\epsilon (d+b)_x)_x,  \label{SWEpos1}\\
  P_t + \left(\frac{P^2}{d}\right)_x + gd(d+b-H)_x &=(v\epsilon (d+b)_x)_x+(\epsilon d v_x)_x.\label{SWEpos2}
\end{align}
Contracting the equations with the entropy variables leads to,
\begin{align}
  U(\uu)_t +F(\uu)_x&=\nonumber\\
  (g(d+b)\epsilon (d+b)_x)_x -g(d+b)_x\epsilon (d+b)_x -(\frac{v^2}{2}\epsilon(d+b)_x)_x + \epsilon vv_x(d+b)_x\nonumber\\
+  (\epsilon v d v_x)_x -\epsilon d v_x^2+(v^2\epsilon (d+b)_x)_x-\epsilon v_xv(d+b)_x&.\nonumber
\end{align}
Clearly, the right-hand side imply a damping. Furthermore, if $v=0$ and $(d+b)=constant$, i.e., a water basin at rest, the artificial diffusion terms in  (\ref{SWEpos1})-(\ref{SWEpos2}) vanish and thereby preserve the steady state, i.e., it is well balanced.

\section{Linear stability}\label{sec:lin_stab}

The non-linear estimate is no guarantee that solutions are not unstable locally around certain flow states (see \cite{GassnerSvard22}). Therefore, we will linearise the system (\ref{eqShallow1})-(\ref{eqShallow2}) and study its stability properties. To this end, we introduce a smooth solution $\tilde v,\tilde d$ and smooth perturbations $v',d'$. (This analysis concerns local stability of time-dependent smooth solutions, not just the steady state.) We let $v=\tilde v+v',d=\tilde d+d'$ in (\ref{eqShallow1})-(\ref{eqShallow2}). The hyperbolic part becomes,
\begin{align*}
  (\tilde d+d')_t+((\tilde d+d')(\tilde v+v'))_x&= 0,\\
  ((\tilde d+d')(\tilde v+v'))_t + ((\tilde d+d')(\tilde v+v')^2)_x + g(\tilde d+d')(\tilde d+d'+b-H)_x &=0.
\end{align*}
Using that $\tilde d,\tilde v$ solves the original system and that $H$ is constant, we obtain
\begin{align*}
  d'_t+(\tilde d v'+\tilde v d' +d'v' )_x&= 0,\\
  (\tilde d v'+\tilde v d' +d'v')_t +(d'(\tilde v^2+2v'\tilde v+(v')^2)
  +\tilde d(2v'\tilde v+(v')^2))_x \\+g\tilde d(d'+b)_x
 +gd'(\tilde d+d'+b)_x &=0.
\end{align*}
Assuming that quadratic perturbations are insignificant compared to first order perturbations gives
\begin{align*}
  d'_t+(\tilde d v'+\tilde v d'  )_x&= 0,\\
  (\tilde d v'+\tilde v d' )_t +(\tilde v^2d'
  +2\tilde d\tilde v v')_x +g\tilde d(d'+b)_x
 +gd'(\tilde d+b)_x &=0.
\end{align*}
Since $\tilde v,\tilde d$ and $b$ are known and bounded functions, we may reduce the system further by ignoring zeroth-order terms in $d',v'$. (Such terms may cause a growth, albeit not unbounded. See \cite{GustafssonKreissOliger}.) This leaves us with,
\begin{align*}
  d'_t+\tilde d v'_x+\tilde v d'_x&= 0,\\
  \tilde d v_t'+\tilde v d'_t  +\tilde v^2d'_x
  +2\tilde d\tilde v v'_x +g\tilde dd'_x&=0,
\end{align*}
leading to,
\begin{align}
  d'_t+\tilde d v'_x+\tilde v d'_x&= 0, \label{lin1}\\
  v_t'-\frac{\tilde v}{\tilde d} (\tilde d v'_x+\tilde v d'_x)  +\frac{\tilde v^2}{\tilde d}d'_x
  +2\tilde v v'_x +gd'_x&=0.\nonumber
\end{align}
We derive an energy estimate for (\ref{lin1}). To reduce notation, we carry it out for the ``frozen coefficient'' case  (see \cite{GustafssonKreissOliger}). That is, we consider fixed constant values of $\tilde d, \tilde v$ and remove the apostrophes. Then (\ref{lin1}) turns into,
\begin{align*}
  d_t+\tilde d v_x+\tilde v d_x&= 0,\\
  v_t+\tilde v v'_x +gd_x&=0.
\end{align*}
We multiply the first equation by $\frac{g}{\tilde d}d$, the second by $v$, sum the resulting equations, and integate in space.
\begin{align*}
  0=\frac{1}{2}(\|\sqrt{\frac{g}{\tilde d}}d\|_2^2+\|v\|_2^2)_t+\int_0^1g dv_x+\frac{g}{\tilde d} \tilde v dd_x +\tilde v vv'_x +gvd_x\,\,dx.
\end{align*}
After integration in space, we are left with $\frac{1}{2}(\|\sqrt{\frac{g}{\tilde d}}d\|^2+\|v\|^2)_t=0$ and we conclude that the shallow-water part is bounded in $L^2$.

Next, we turn to the dispersive terms. The third derivative term is linearised as,
\begin{align*}
  (\gamma v_{xx})_x&=(\gamma (\tilde v+v')_{xx})_x=(\gamma \tilde v_{xx})_x+(\gamma v'_{xx})_x,\\
(\gamma v_{x})_{xx}&=(\gamma (\tilde v+v')_{x})_{xx}=(\gamma \tilde v_{x})_{xx}+(\gamma v'_{x})_{xx}.
\end{align*}
We neglect the forcing terms that do not affect well-posedness. Since the other dispersive terms are also linear, we obtain the following adjustment of the previous estimate,
\begin{align*}
  0=\frac{1}{2}(\|\sqrt{\frac{g}{\tilde d}}d\|^2+\|v\|^2)_t-\int_0^1\left(v(\gamma v_{xx})_x+v(\gamma v_x)_{xx}+v(\beta v_x)_{xt}\right) \,\,dx.
\end{align*}
(We have dropped the apostrophes.) It follows that $\frac{1}{2}(\|\sqrt{\frac{g}{\tilde d}}d\|^2+\|v\|^2)_t+\beta\left(\frac{1}{2} (v_x)^2\right)_t=0$.

Finally, we turn to the artificial diffusion. The term in the $d$-equation is : $\epsilon(\tilde d+d'+b)_x)_x$. Only the $d'$ must be kept in the linear stability analysis. In the momentum equation, the diffusion is non-linear:
\begin{align*}
  (v\epsilon (d+b)_x)_x+(\delta v_x)_x=
  ((\tilde v+v')\epsilon (\tilde d+d'+b)_x)_x+(\delta (\tilde v+v')_x).
\end{align*}
Ignoring forcing terms, we are left with,
\begin{align*}
  (v'\epsilon (\tilde d+b)_x)_x+  (\tilde v\epsilon d'_x)_x+(\delta v'_x)_x.
\end{align*}
These terms augment the energy estimate as
\begin{align}
  \int_0^1 \frac{g}{\tilde d}d'\epsilon(\tilde d'_x)_x+  v'(v'\epsilon (\tilde d+b)_x)_x+  v'(\tilde v\epsilon d'_x)_x+v'(\delta v'_x)_x\,\,dx&=\\
  \int_0^1 -\epsilon \frac{g}{\tilde d}  \frac{(d'_x)^2}{2}  -v'_xv'\epsilon (\tilde d+b)_x-  v'_x(\tilde v\epsilon d'_x)-\delta (v'_x)^2\,\,dx&.
\end{align}
The second term on the last row is linearly bounded since
\begin{align*}
\int -v'_xv'\epsilon (\tilde d+b)_x\,dx= \int \frac{(v')^2}{2}\epsilon (\tilde d+b)_{xx},
\end{align*}
which only contributes with at most an exponential growth. 

The third term can be bounded by the first and the last by observing that
\begin{align}
  \int_0^1 -  v'_x(\tilde v\epsilon d'_x)\,\,dx&\leq \int_0^1 \frac{1}{2}(\sqrt{\epsilon}v'_x\tilde v\sqrt{\frac{\tilde d}{g}})^2+ \frac{1}{2}(\sqrt{\epsilon}\sqrt{\frac{g}{\tilde d}}d'_x)^2.
\end{align}
Hence, we must require  $\delta\geq \frac{1}{2}\frac{\tilde d}{g}\epsilon |\tilde v|$ for linear stability. 

These brief considerations demonstrate that the problem is linearly stable in addition to satisfying the non-linear entropy estimate. However, it should be noted that the linear energy estimate does not imply that the energy is always decaying. The zeroth-order terms and forcing functions that we omitted may induce a considerable local growth of the solution. Likewise, non-constant $\hat \alpha,\hat \beta,\hat\gamma$ may induce a linear growth. However, in all these case, the growth remains within the limits of linear stability.

\section{Numerical scheme}

A numerical scheme for (\ref{eqShallow1})-(\ref{eqShallow2}) should satisfy the discrete counterpart of a
nonlinear entropy estimate in order to retain the same stability properties as the PDE itself
(such schemes are termed \emph{entropy stable}). Furthermore, it should be well balanced,
meaning that a ``lake at rest'' is a steady state solution to the numerical scheme.
We will also propose entropy-stable artificial-diffusion terms that can be switched on where necessary.

For the shallow-water equations, well balanced and entropy-conservative schemes  have already been derived
(see \cite{Fjordholm_etal11} and \cite{Wintermeyer_etal17}, Theorem 1).
We begin by verifying their results in a 1-D finite volume setting. To this end, we need the notation:
\begin{align}
  \bar a_{i+1/2}&=\frac{a_{i+1}+a_i}{2}, \nonumber\\
  \overline{a^2}_{i+1/2}&=\frac{a^2_{i+1}+a^2_i}{2},\label{num_rels}\\
  \Delta_+ a_{i}&=\Delta_-a_{i+1}=a_{i+1}-a_i,\nonumber\\
  D_-a_{i+1}&=h^{-1}\Delta_-a_{i+1}=D_+a_i.\nonumber
\end{align}
For simplicity, we assume a constant grid spacing with step size $h$ and a periodic domain. The primal grid is thus $x_i=ih$ and the dual grid $x_{i+1/2}=x_i+h/2$. We consider the finite volume scheme
\begin{align}
(d_i)_t+\frac{P_{i+1/2}-P_{i-1/2}}{h}&=0,\label{SWscheme1}\\
(P_i)_t+ \frac{Q_{i+1/2}-Q_{i-1/2}}{h}+g\frac{1}{2}\left(\bar d_{i+1/2}D_+b_i+\bar d_{i-1/2}D_-b_i\right)&=0,\label{SWscheme2}
\end{align}
where
\begin{align}
  P_{i+1/2}&=\bar d_{i+1/2}\bar v_{i+1/2},\label{num_flux}\\
  Q_{i+1/2}&=\bar d_{i+1/2}(\bar v_{i+1/2})^2+\frac{1}{2}g\overline{d^2}_{i+1/2}\nonumber.
\end{align}
We also need the entropy potential, that is given by
\begin{align}
\Psi=\frac{1}{2g}(w_1^2w_2+w_1w_2^3)+\frac{1}{8g}w_2^5=\frac{1}{2}gd^2v,\nonumber
\end{align}
where $w_j$, $j=1,2$ denote the components of the entropy variables (\ref{ent_var}).

The scheme is {\bf well-balanced}, if a basin of water at rest, with varying bathymetry, is a steady state solution of the scheme.  That is, $d+b=constant$ and $v=0$ should be a solution of the scheme. Inserting  $v=0$ in (\ref{SWscheme1}) yields $P_{i+1/2}=0$, which in turn implies that $(d_i)_t=0$.
Similarly, having $v=0$ results in
\begin{align*}
Q_{i+1/2}=\frac{1}{2}g\overline{d^2}_{i+1/2},
\end{align*}
which inserted in (\ref{SWscheme2}) yields,
\begin{align*}
  (P_i)_t +\frac{1}{2}g\frac{\overline{d^2}_{i+1/2}-\overline{d^2}_{i-1/2}}{h}+g\frac{1}{2}\left(\bar d_{i+1/2}D_+b_i+\bar d_{i-1/2}D_-b_i\right)&=0.
\end{align*}
Using (\ref{num_rels}), we obtain
\begin{align*}  
  (P_i)_t +\frac{1}{2}g\frac{d^2_{i+1}+d_i^2-d_i^2-d^2_{i-1}}{2h}+g\frac{1}{2}\left(\bar d_{i+1/2}D_+b_i+\bar d_{i-1/2}D_-b_i\right)&=0,
\end{align*}
and with $d_i=C-b_i$, where $C$ is a constant, we have
  \begin{align*}
  (P_i)_t +\frac{1}{2}g\left(\frac{d_{i+1}+d_{i}}{2}\frac{d_{i+1}-d_{i}}{h}+\frac{d_{i+1}+d_{i}}{2}D_+(C-d_i)\right)&\\
    +\frac{1}{2}g\left(\frac{d_{i}+d_{i-1}}{2}\frac{d_{i}-d_{i-1}}{h}+g\frac{d_{i}+d_{i-1}}{2}D_-(C-d_i)\right)  &=0.
\end{align*}
We arrive at  $(P_i)_t=0$. Thus, the basin at rest is a steady state solution and the scheme is thus well balanced.

For the entropy analysis, we need the standard relations:
\begin{align}
  \overline{dv}_{1/2}&=\bar v_{1/2}\bar d_{1/2}+\frac{1}{4}\Delta_+v_0  \Delta_+ d_0,\nonumber\\
  D_-(v_1d_1)&=\bar v_{1/2}D_+d_i+\bar d_{1/2}D_+v_i,\label{standard_rels}\\
  D_-(d_{i+1}v_{i+1})&=d_{i+1}D_-v_{i+1}+v_iD_-d_{i+1}.\nonumber
\end{align}

Turning to {\bf entropy conservation}, we introduce the discrete entropy function $U_i=U(\uu_i)=(\frac{P^2_i}{2d_i^2}+\frac{gd_i^2}{2}+gd_ib_i)$ where $\uu_i=(d_i, P_i)^T$. Then the discrete entropy variables are, $\partial_{\uu_i}U_i=\ww^T_i=\left(g(d_i+b_i)-\frac{P^2_i}{2d^2_i}, \frac{P_i}{d_i}\right)$.

In analogy with the continuous case, we contract the scheme (\ref{SWscheme1})-(\ref{SWscheme2}) with $h\ww^T_i$,
\begin{align}
h (w_1)_i\left((d_i)_t+\frac{P_{i+1/2}-P_{i-1/2}}{h}\right)\nonumber\\
 + h(w_2)_i\left((P_i)_t+ \frac{Q_{i+1/2}-Q_{i-1/2}}{h}+g\frac{1}{2}\left(\bar d_{i+1/2}D_+b_i+\bar d_{i-1/2}D_-b_i\right)\right)&=0.\label{ent_est}
\end{align}
The time derivative terms are combined as in the continuous case to obtain
\begin{align*}
  h(U_i)_t+(w_1)_i(P_{i+1/2}-P_{i+1/2}) \\
   + (w_2)_i\left((Q_{i+1/2}-Q_{i+1/2})+g\frac{1}{2}\left(\bar d_{i+1/2}\Delta_+b_i+\bar d_{i-1/2}\Delta_-b_i\right)\right)&=0.
\end{align*}
We recast the equation in the standard way (\cite{Tadmor03}),
\begin{align*}
  h(U_i)_t+\overline{w_1}_{i+1/2}P_{i+1/2}-\frac{1}{2}\Delta_+{w_1}_{i}P_{i+1/2}
  -\overline{w_1}_{i-1/2}P_{i-1/2}-\frac{1}{2}\Delta_-{w_1}_{i}P_{i-1/2}  \\
  +\overline{w_2}_{i+1/2}Q_{i+1/2}-\frac{1}{2}\Delta_+{w_2}_{i}Q_{i+1/2}
  -\overline{w_2}_{i-1/2}Q_{i+1/2}-\frac{1}{2}\Delta_-{w_2}_{i}Q_{i-1/2}\\
  + \frac{g}{2}\overline{w_2}_{i+1/2}\bar d_{i+1/2}\Delta_+b_i
  - \frac{g}{2}\Delta_+{w_2}_{i}\bar d_{i+1/2}\Delta_+b_i
  + \frac{g}{2}\overline{w_2}_{i-1/2}\bar d_{i-1/2}\Delta_-b_i
  + \frac{g}{2}\Delta_-{w_2}_{i}\bar d_{i-1/2}\Delta_-b_i
&=0.
\end{align*}

By (\ref{ent_var}), $w_2=v$ and we define the entropy flux as,
\begin{align}
F_{i+1/2}=\overline{( w_1)}_{i+1/2} P_{i+1/2}+\overline{( w_2)}_{i+1/2} Q_{i+1/2}- \frac{g}{2}\Delta_+v_i\bar d_{i+1/2}\Delta_+b_i- \bar\Psi_{i+1/2},\label{numFflux}
\end{align}
where $\Psi_i=\frac{g}{2}P_id_i$ is the entropy potential. (It is straightforward to verify that the numerical entropy flux (\ref{numFflux}) is consistent with (\ref{F}).)
This results in,
\begin{align}
  h (U_i)_t+F_{i+1/2}-F_{i+1/2}\nonumber\\
  -\frac{1}{2}\Delta_+ ( w_1)_{i}P_{i+1/2}
-\frac{1}{2}\Delta_+ (w_2)_{i}Q_{i+1/2} +\Delta_+ \Psi_{i}\nonumber\\
  -\frac{1}{2}\Delta_- (w_1)_{i}P_{i-1/2}
  -\frac{1}{2}\Delta_- (w_1)_{i}Q_{i+1/2}+\Delta_- \Psi_{i}\nonumber\\
    + \frac{g}{2}\overline{w_2}_{i+1/2}\bar d_{i+1/2}\Delta_+b_i
  + \frac{g}{2}\overline{w_2}_{i-1/2}\bar d_{i-1/2}\Delta_-b_i&=0.\nonumber
\end{align}
For entropy conservation, we must require that
\begin{align}
\frac{1}{2}\Delta_+ ( w_1)_{i}P_{i+1/2}+
\frac{1}{2}\Delta_+ ( w_2)_{i}Q_{i+1/2} - \frac{g}{2}\overline{w_2}_{i+1/2}\bar d_{i+1/2}\Delta_+b_i=\frac{1}{2}\Delta \Psi_{i+1/2}. \label{shuffle}
\end{align}
To verify that this holds, insert the explicit form of the numerical fluxes (\ref{num_flux}) and entropy variables (\ref{ent_var}) in the left-hand side of (\ref{shuffle}):
\begin{align*}
A=  \frac{1}{2}\Delta_+ ( w_1)_{i}P_{i+1/2} + \frac{1}{2}\Delta_+ ( w_2)_{i}Q_{i+1/2} - \frac{g}{2}\overline{v}_{i+1/2}\bar d_{i+1/2}\Delta_+b_i&=\\
  \frac{1}{2}\left(g(\Delta_+ d_{i}+\Delta_+ b_{i})-\Delta_+\left(\frac{v^2}{2}\right)_{i}\right)\bar d_{i+1/2}\bar v_{i+1/2}\\
  +\frac{1}{2}\Delta_+  v_{i}\left(\bar d_{i+1/2}(\bar v_{i+1/2})^2+\frac{1}{2}g\overline{ d^2}_{i+1/2}\right) -\frac{g}{2}\overline{v}_{i+1/2}\bar d_{i+1/2}\Delta_+b_i.
\end{align*}
We simplify to obtain,
\begin{align*}
2A=  
  (g(\Delta_+ d_{i}+\Delta_+ b_{i})-\bar v_{i+1/2}\Delta_+ v_{i})\bar d_{i+1/2}\bar v_{i+1/2}\\
+\Delta_+  v_{i}(\bar d_{i+1/2}(\bar v_{i+1/2})^2+\frac{1}{2}g\overline{ d^2}_{i+1/2})
- g\overline{v}_{i+1/2}\bar d_{i+1/2}\Delta_+b_i
&=\\
  (g(\Delta_+ d_{i}+\Delta_+ b_{i}))\bar d_{i+1/2}\bar v_{i+1/2}
+\Delta_+  v_{i}\frac{1}{2}g\overline{ d^2}_{i+1/2}
- g\overline{v}_{i+1/2}\bar d_{i+1/2}\Delta_+b_i
&=\\
  \frac{1}{2}g(\Delta_+ d^2_{i})\bar v_{i+1/2}
  +\Delta_+  v_{i}\frac{1}{2}g\overline{ d^2}_{i+1/2}
  &=\\
  =\frac{1}{2}g\Delta_+ (vd^2)_{i}& =
  \Delta_+\Psi_i.
\end{align*}
Hence, (\ref{shuffle}) holds. Collecting the results, we have recast (\ref{ent_est}) into,
\begin{align}
  h (U_i)_t+F_{i+1/2}-F_{i+1/2}&=0,\nonumber
\end{align}
and verified the result in \cite{Fjordholm_etal11}.

\vspace{0.25cm}

In accordance with the analysis in Section \ref{sec:art_diff}, we consider artificial diffusion of the following form,
\begin{align}
(d_i)_t+\frac{P_{i+1/2}-P_{i-1/2}}{h}=&D_-\lambda_{i+1/2}\Delta_+(d_i+b_i)\label{SW_AD1},\\
  (P_i)_t+ \frac{Q_{i+1/2}-Q_{i-1/2}}{h}+g\frac{1}{2}\left(\bar d_{i+1/2}D_+b_i+\bar d_{i-1/2}D_-b_i\right)=&D_-\lambda_{i+1/2}v_{i+1/2}\Delta_+(d_i+b_i)\label{SW_AD2}\\
  &+D_-\lambda_{i+1/2}d_{i+1/2}\Delta_+v_i.\nonumber
\end{align}
Here, we choose the same diffusive coefficient, $\lambda_{i+1/2}\geq 0$, for both terms in the momentum equations since then it collapses to the naive momentum diffusion when $b=constant$. The well-balancedness of the artificial diffusion is trivially verified by inserting $d_i+b_i=constant$ and $v_i=0$.

Turning to entropy stability, we contract the scheme with the entropy variables to obtain
\begin{align}
  h (U_i)_t+\tilde F_{i+1/2}-\tilde F_{i+1/2} &=AD,\nonumber
\end{align}
where
\begin{align*}
\tilde F_{i+1/2}=& F_{i+1/2}-\overline{( w_1)}_{i+1/2} \lambda_{i+1/2}\Delta_+(d_i+b_i)\\
 &- \overline{( w_2)}_{i+1/2} \left(\lambda_{i+1/2}v_{i+1/2}\Delta_+(d_i+b_i)
  +\lambda_{i+1/2}d_{i+1/2}\Delta_+v_i \right),
\end{align*}
and $AD=AD_{i+1/2}+AD_{i-1/2}$ is the remainder from the artificial diffusion terms. The explicit form of these terms is given by,
\begin{align*}
AD_{i+1/2}=&-\left(g\Delta_+ (d_{i}+ b_{i})-\bar v_{i+1/2}\Delta_+ v_{i}\right)\lambda_{i+1/2}\Delta_+(d_i+b_i)\\
&-\Delta_+v_i \left(\lambda_{i+1/2}v_{i+1/2}\Delta_+(d_i+b_i)
+\lambda_{i+1/2}d_{i+1/2}\Delta_+v_i \right)\\
=&
-g\lambda_{i+1/2}(\Delta_+(d_i+b_i))^2-\lambda_{i+1/2}d_{i+1/2}(\Delta_+v_i)^2\leq 0.
\end{align*}
Hence, the scheme is dissipating entropy.

\vspace{0.25cm}


By design Boussinesq systems are supposed to be purely dispersive.  However, there are some issues with entropy conservative schemes of the type proposed above: They are not necessarily locally linearly stable. (See \cite{GassnerSvard22}.) This can easily be seen by considering the scheme for the continuity equation,
\begin{align}
(d_i)_t+\frac{\bar d_{i+1/2}\bar v_{i+1/2}-\bar d_{i-1/2}\bar v_{i-1/2}}{h}=&D_-\lambda_{i+1/2}\Delta_+(d_i+b_i)\label{cont}.
\end{align}
Assuming a bounded and smooth velocity, this is a linear advection equation. The marginally linearly stable approximation for such an equation is:
\begin{align*}
(d_i)_t+\frac{d_{i+1}v_{i+1}- d_{i-1} v_{i-1}}{2h}=&0.
\end{align*}
It is straightforward to recast (\ref{cont}) into this form,
\begin{align}
  (d_i)_t+\frac{d_{i+1} v_{i+1}-d_{i-1} v_{i-1}}{2h}
-\frac{1}{4}D_-(\Delta_+v_i\Delta_+d_i)
  =&D_-\lambda_{i+1/2}\Delta_+(d_i+b_i).\label{cont}
\end{align}
Since $\Delta_+v_i$ can take any sign, the scheme may be locally anti-diffusive, which in turn is linearly unstable. By choosing $\lambda_{i+1/2}\geq \frac{1}{4}|\Delta_+v_i|$, these instabilities are suppressed. (Note that $\Delta_+v_i$ is $\mathcal{O}(h)$ for smooth solutions such that the scheme is still second-order accurate with this choice of $\lambda_{i+1/2}$.) 
However, to ensure grid convergence, we typically need to add a first-order artificial diffusion due to the non-linear structure of the equations.

\subsection{The Boussinesq terms}

Let $D_2=D_+D_-$ and note that $D_2=D_-D_+$.
Furthermore, $D_0=\frac{1}{2}(D_++D_-)$ is the standard central difference stencil.

We add the Boussinesq terms to the shallow-water scheme (\ref{SW_AD1})-(\ref{SW_AD2}):
\begin{align}
  (d_i)_t+D_-\tilde P_{i+1/2}&=
D_3^d(\hat\alpha,d,b),\label{B1}
  \\
  (P_i)_t+ D_-\tilde Q_{i+1/2} +g\frac{1}{2}\left(\bar d_{i+1/2}D_+b_i+\bar d_{i-1/2}D_-b_i\right)&=\label{B2}\\
D_3^P(\hat\alpha,v,d,b)
+\frac{1}{2}D_-(\hat\beta_i D_+v_i)_{t}+ \frac{1}{2}D_+(\hat\beta_i D_-v_i)_{t}\nonumber\\ 
+ \frac{1}{2}D_+(\hat\gamma_iD_-D_+v_i)+ \frac{1}{2}D_-D_+(\hat\gamma_iD_-v_i)\nonumber,
\end{align}
where
\begin{align*}
\tilde P_{i+1/2}&=P_{i+1/2}-\lambda_{i+1/2}\Delta_+(d_i+b_i)),\\
\tilde Q_{i+1/2}&=Q_{i+1/2}-\lambda_{i+1/2}v_{i+1/2}\Delta_+(d_i+b_i)-\lambda_{i+1/2}d_{i+1/2}\Delta_+v_i,
\end{align*}
and $\{\hat \alpha,\hat\beta,\hat\gamma\}_i$ are the obvious pointwise projections of $\{\hat \alpha,\hat\beta,\hat\gamma\}$ at $x_i$. Moreover,
\begin{align*}
  D_3^d(\hat\alpha,d,b)=&\frac{1}{2}\left(D_-(\hat \alpha_{i} D_+(\hat \alpha_{i} D_-(d+b)_i))+D_+(\hat \alpha_{i} D_-(\hat \alpha_{i} D_+(d+b)_i))\right),\\
  D_3^P(\hat\alpha,v,d,b)=&\frac{1}{2}\left(D_- (\bar v_{i-1/2}\hat \alpha_{i} D_+(\hat \alpha_{i} D_-(d+b)_i))\right.\\
& \left. +D_+(\bar v_{i+1/2}\hat \alpha_{i} D_-(\hat \alpha_{i} D_+(d+b)_i))\right).
\end{align*}

It is evident that the additional terms are \emph{well-balanced} and we proceed to demonstrate entropy boundedness. We already know that the left-hand side is an entropy-stable discretisation and we focus on the Boussinesq terms. As before, we contract with the entropy variables and we consider the $\hat \alpha,\hat\beta, \hat \gamma$ terms separately. (Denoted $T^\alpha_i,T^\gamma_i,T^\beta_i$.)

We begin with the easiest one:
\begin{align*}
  2T^\beta_i=v_iD_-(\hat\beta_i D_+v_i)_{t}+v_iD_+(\hat\beta_i D_-v_i)_{t}&=\\
  D_-(v_i\hat \beta_iD_+(v_i)_{t})-(D_+v_{i-1})\hat \beta_{i-1}(D_+ v_{i-1})_{t}&\\
+ D_+(v_i\hat \beta_iD_-(v_i)_{t})-(D_-v_{i+1})\hat \beta_{i+1}(D_- v_{i+1})_{t}  &=\\
  D_-(v_i\hat \beta_iD_+(v_i)_{t})-((\sqrt{\hat \beta_{i-1}}D_+ v_{i-1})^2)_{t}&\\
+ D_+(v_i\hat \beta_iD_-(v_i)_{t})-((\sqrt{\hat \beta_{i+1}}D_- v_{i+1})^2)_{t}  &.
\end{align*}
Since we are only interested in $\sum_{i=0}^NhT^{\gamma}_i$ on a periodic domain, we simplify the calculations by using a standard summation-by-parts rules (see e.g. \cite{GustafssonKreissOliger}). We obtain,
\begin{align}
  \sum_i h 2T^\gamma_i=\sum_i v_i\left(D_+(\hat\gamma_iD_-D_+v_i)+D_-D_+(\hat\gamma_iD_-v_i)\right)&=\nonumber \\
  \sum_i -(D_-v_i)\hat\gamma_i(D_-D_+v_i)+(D_-D_+v_i)(\hat\gamma_iD_-v_i)&=0\nonumber .
\end{align}
We omit the proof for $T_{\alpha}$.
  Upon summation in space and integration in time an estimate of the form
\begin{align*}
  \sum_{i=1}^N hU_i|_{t=T}+(\sqrt{\hat \beta_{i+1}}D_- v_{i+1})^2|_{t=T}\leq Constant,
\end{align*}
  is obtained for any bounded time $T$. Thus, the scheme is entropy stable.

\subsection{  Regularity of the bathymetry function }

In the non-linear entropy analysis (section \ref{sec:var_bat}), it is not necessary to assume any regularity restrictions on $b(x)$. (Recall that $b(x)$ enters the dispersive coefficients via the still water depth: $h(x)=H-b(x)$.)

However, in the linear stability analysis in section \ref{sec:lin_stab}, one has to assume that the dispersive coefficients, i.e., $b(x)$, have at least two bounded derivatives. (This ensures that some of the neglected terms are bounded.) Furthermore, we emphasise again that these forcing terms will contribute to a local exponential (but not unstable \cite{GustafssonKreissOliger}) growth. Thus the system will be more sensitive to numerical errors in the vicinity of a strongly varying bathymetry.

Although, we do not have a complete non-linear well-posedness analysis at hand, it seems likely that some regularity of $b(x)$ will be required to ensure numerical grid convergence. This is easily achieved by mollifying $b(x)$. This is hardly a very restrictive requirement as the bathymetry is presumable only approximately known anyway.

A simple, yet effective, way to regularise $b(x)$ mollify $b(x)$ with the indicator function $\phi(x,y)$. The indicator function $\phi(x,y)$ is 0 everywhere but $x-\delta\leq y\leq x+\delta $. If $b(x)$ is piecewise continuous, then  $b_0(x)=\int_\Omega \phi(x,y)b(y)dx$ is continous. Repeated application of the mollifier increases the regularity by one each time.

To demonstrate the robustness of the system, we do not mollify the bathymetry in the numerical examples below. (However, we have verified that the mollificiation procedure above gives very similar results.)

\section{Time stepping scheme}

The semi-discrete system (\ref{B1}) and (\ref{B2}) can in principle be advanced explicitly in time but the third derivatives constrain the time-step size severly making this approach unfeasible for practical problems.  A fully implicit scheme allows large time steps but is complicated to code and computationally expensive to solve due to the non-linear convective terms. However, a fully implicit time stepping scheme is unnecessary since we anyway need to resolve the convective time scale to obtain accurate solutions. Hence, an implicit-explicit (IMEX) scheme is preferable. Since the focus of this paper is \emph{not} on time stepping schemes, we will only propose a temporal scheme that is formally first-order accurate in time. We will also restrict ourselves to systems with $\alpha=0$ when deriving the IMEX scheme. We emphasise that we have verified that the general scheme  (\ref{B1}) and (\ref{B2}) is stable (with $\alpha\neq 0$) using explicit time stepping.)

We discretise time as $t^n=nk$, where $k$ is the time step. (The generalisation to variable time steps is trivial.) Furthermore, the approximations at time $t^n$  are denoted with a superscript. For instance, $P^n_i$ is the approximation of the momentum at $(x_i,t^n)$. We introduce the time difference operator,
\begin{align*}
D_T\xi_i^n=\frac{\xi_i^{n+1}-\xi^n_i}{k},
\end{align*}
where $\xi$ represents a variable or flux. In equation, (\ref{B2}), there are time derivatives in both $P$ and $v$. To obtain a solvable scheme, we make the following reformulation:
\begin{align*}
  P_t=(dv)_t = d_tv + dv_t.
\end{align*}
Thus, a fully discrete scheme for (\ref{B1})-(\ref{B2}), with $\alpha=0$, is
\begin{align}
  D_Td_i^n+D_-\tilde P^n_{i+1/2}&=0\label{F1}\\
v_i^nD_Td_i^n + d_i^nD_Tv_i^n+ D_-\tilde Q^n_{i+1/2} +g\frac{1}{2}\left(\bar d^n_{i+1/2}D_+b_i+\bar d^n_{i-1/2}D_-b_i\right)&=\label{F2}\\
+\frac{1}{2}D_TD_-(\hat\beta_i D_+v_i^n)+ \frac{1}{2}D_TD_+(\hat\beta_i D_-v^n_i)\nonumber\\ 
+ \frac{1}{2}D_+(\hat\gamma_iD_-D_+v^{n+1}_i)+ \frac{1}{2}D_-D_+(\hat\gamma_iD_-v^{n+1}_i)\nonumber.
\end{align}
Equation (\ref{F1}) is a standard explicit Euler scheme. In (\ref{F2}), the nonlinear shallow-water part, is taken explicitly and the third-derivative  $\gamma$-terms are taken implicitly. Using (\ref{F1}), we recast (\ref{F2}) to obtain the scheme
\begin{align}
   D_Td_i^n+D_-\tilde P^n_{i+1/2}&=0\label{FU1}\\
   d_i^nD_Tv_i^n -\frac{1}{2}D_TD_-(\hat\beta_i D_+v_i^n)- \frac{1}{2}D_TD_+(\hat\beta_i D_-v^n_i),\nonumber\\
- \frac{1}{2}D_+(\hat\gamma_iD_-D_+v^{n+1}_i)- \frac{1}{2}D_-D_+(\hat\gamma_iD_-v^{n+1}_i)   &=\label{FU2}\\
 v_i^nD_-\tilde P^n_{i+1/2}   - D_-\tilde Q^n_{i+1/2} -g\frac{1}{2}\left(\bar d^n_{i+1/2}D_+b_i+\bar d^n_{i-1/2}D_-b_i\right).
\nonumber
\end{align}
This scheme is straightforward to implement:
\begin{enumerate}
\item Take a step with (\ref{FU1}) and save $D_-\tilde P^n_{i+1/2}$.
\item Using $D_-\tilde P^n_{i+1/2}$, the right-hand side of (\ref{FU2}) can be computed explicitly.
  \item The resulting equation is a linear system of equations (across all points, $i$, in space ) and it is solved for $v^{n+1}_i$.
\end{enumerate}

\section{Numerical tests}

The exact scheme (\ref{FU1})-(\ref{FU2}) without any additional numerical dissipation,
has been implemented in Julia \cite{julia}. To demonstrate the robustness of the current approach,
we run all cases without mollifying the bathymetry. That is, corners are sharp and
jumps are true discontinuities. The first problem we consider models waves that run past
a trapezoidal sill. We run the code on a sequence of finer grids to demonstrate grid convergence.
This appproach allows us to single out modelling errors from numerical errors.
Next, we run different dispersive coefficients,
namely both set 3 and set 4, and compare the results with experimental data \cite{Dingemans94}.

To further test the robustness, we run two more problems. The first is a sharp triangular sill. Apart from testing strongly varying bathymetries,  it also demonstrates the automatic reduction of the dispersive terms in shallow water. The second is a bathymetry featuring a cavity where $b(x)$ is discontinuous.

\subsection{Validation}\label{sec:dinge}

To validate the model, we run Dingeman's experimental (\cite{Dingemans94,Dingemans97}) setup for waves passing a trapezoidal bar. In experiment, the wave maker is situated at $x=0$ and the depth is 80 cm. (The depth is actually $86$cm for the first $6$ meters but that is claimed not to affect the results.) The bathymetry is then given by:
\begin{itemize}
\item $x<11.01$ meters; 80 cm constant depth.
\item $11.01\leq x < 23.04$, linear slope to depth 20 cm.
\item  $23.04\leq  x< 27.04$; 20 cm constant depth.
\item $27.04\leq x <33.07$, linear slope to depth 80 cm.
\item $x\geq 33.07$; 80 cm constant depth.
\end{itemize}
The wave maker generated waves with an ampltidude of $2$cm. Measurements of the wave height were taken
at two sets of locations.
\begin{align*}
x^A&=\{3.04,9.44, 20.04, 26.04, 30.44, 37.04\},\\
x^B&= \{7.04,9.44, 24.04, 28.04,  33.64, 41.04\}.
\end{align*}
Herein, we use data obtained at $x_A$ for comparisons.

In our simulations, we do not model a wave maker. Instead, we generate a wave train with eight crests
that have an amplitude of $2$cm. The initial wave is generated using the dispersive relation of the
Euler equations as follows:
\begin{align*}
  d&=A\cos{Kx},\nonumber\\
  v&=\sqrt{ \frac{g}{K}\tanh{K d_0}}\frac{d}{d_0},
\end{align*}
where $d_0=0.8$  is the still water depth. $A=0.02$ is the amplitude. The gravitational constant is set to $g=9.81$. The wave number $K$ is the solution of the dispersive relation,
\begin{align*}
\omega^2=g K \tanh(K d_0)
\end{align*}
where $\omega=\frac{2\pi}{T}$ and $T=2.02\sqrt{2}$. At both ends, we smoothly let the amplitude approach zero. The wave train is placed well in front of the sill such that there are no initial interactions. Hence, we need to extend the computational domain compared to the experimental setup. We use the domain $x\in [-138,46]$ and we make it periodic, since we have not analyzed boundary conditions properly. We divide the computational domain into $N$ equidistant cells.

\begin{figure}[h]
\centering
\includegraphics[width=6cm]{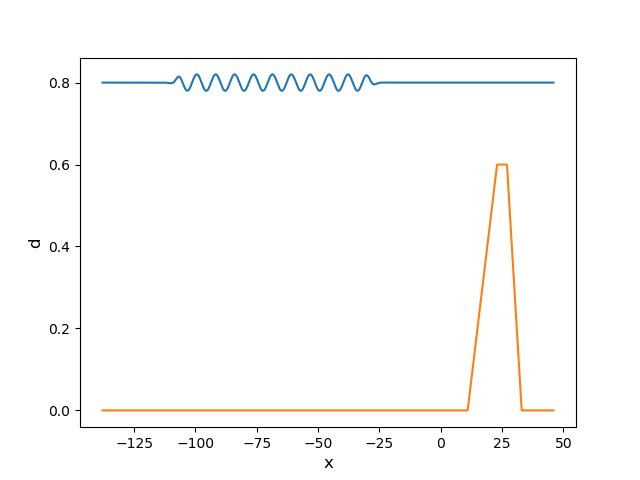}\label{fig:dinge_init}
\caption{initial data and bottom topography}\label{fig:dinge}
\end{figure}

The bottom topography and the initial data (for the depth variable) is shown in Fig. \ref{fig:dinge_init}.

We begin by running the code with the dispersive parameters (\ref{set3}) ``Set 3'' on grids with spacing $h=\{0.1,0.05,0.01,0.005\}$. On the coarsest grid, there are $1840$ points in space and on the finest $36800$. For the initial data, this corresponds to approximately $75$ grid points per wavelength on the coarsest grid and $1500$ points on the finest.
 
All computations were run with constant time steps $k=CFL\cdot h$, where $CFL=0.2$ till $T=70$. The artificial diffusion coefficient is taken to be a constant, $\lambda_{i+1/2}=\lambda$ for all $i$. Since the initial velocity has an amplitude of $\approx  0.0654$ and the oscillations increase somewhat during the interaction with the sill, we take $\lambda=0.1$. (Surely, one can choose $\lambda_{i+1/2}$ in a more sophisticated way, but we leave that as future work.)

Finally, there is a translation between the time variable in the computions and experiments. Hence, we shift our solutions by a constant. In order to see the grid convergence it is chosen (by ocular inspection) for the finest grid. (All simulations presented herein are translated with the same constant in time.)

\begin{figure}[ht]
\centering
\subfigure[$h=0.1$]{\includegraphics[width=6cm]{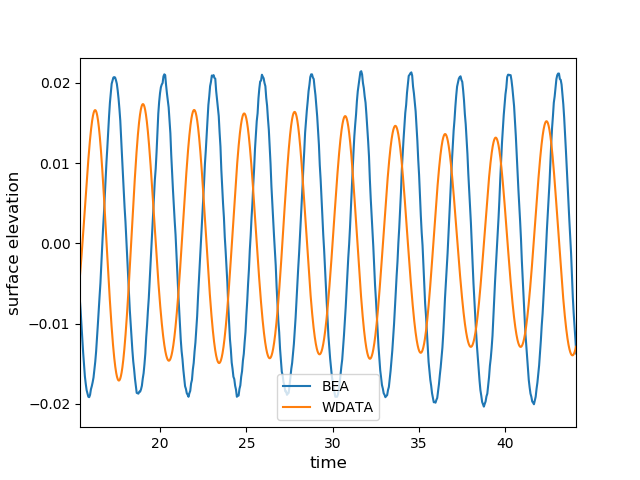}}
\subfigure[$h=0.1$]{\includegraphics[width=6cm]{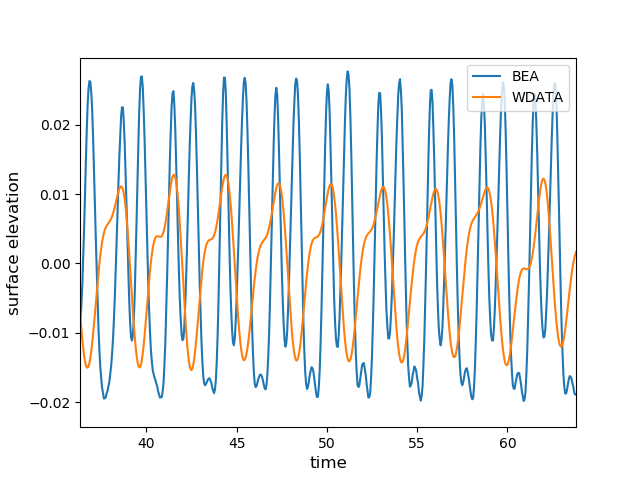}}
\subfigure[$h=0.05$]{\includegraphics[width=6cm]{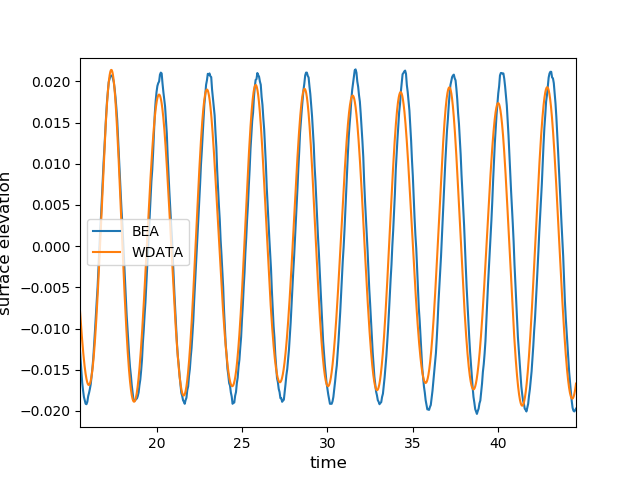}}
\subfigure[$h=0.05$]{\includegraphics[width=6cm]{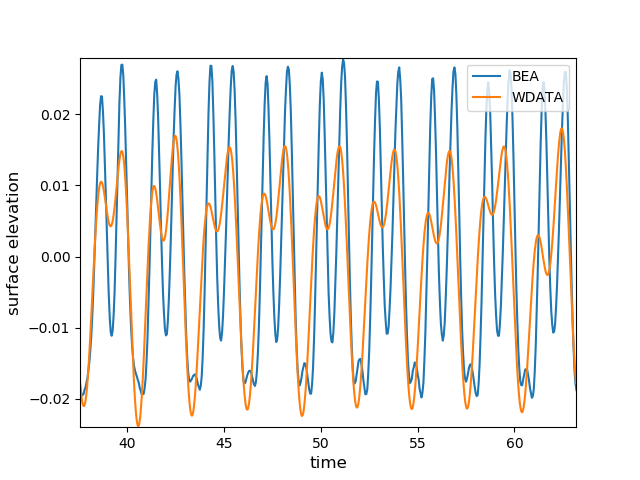}}
\subfigure[$h=0.01$]{\includegraphics[width=6cm]{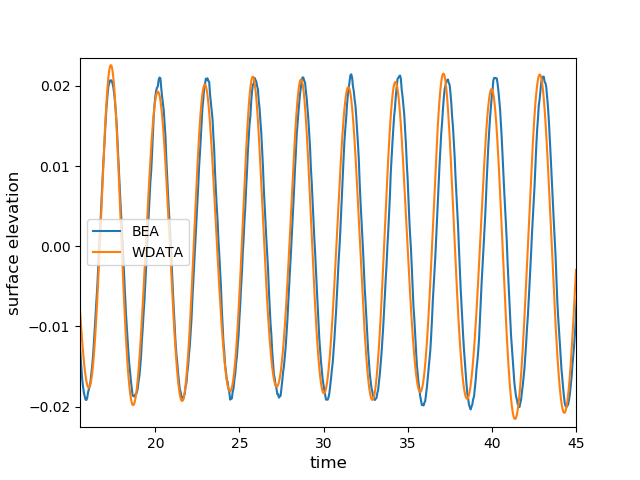}}
\subfigure[$h=0.01$]{\includegraphics[width=6cm]{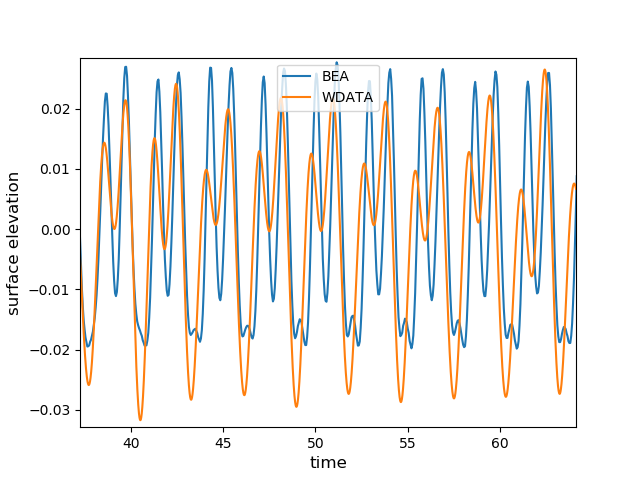}}
\subfigure[$h=0.005$]{\includegraphics[width=6cm]{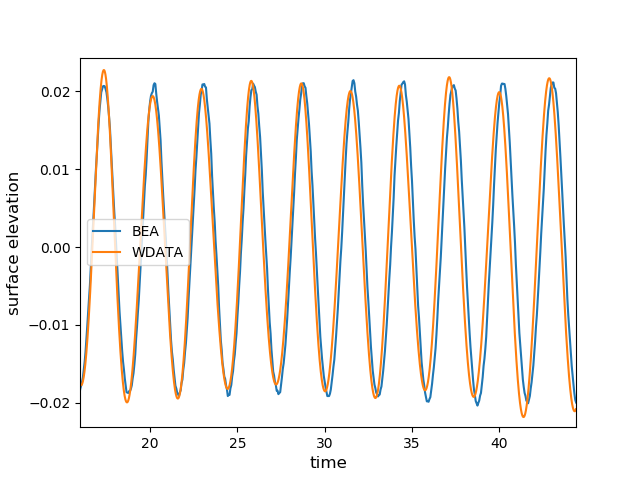}}
\subfigure[$h=0.005$]{\includegraphics[width=6cm]{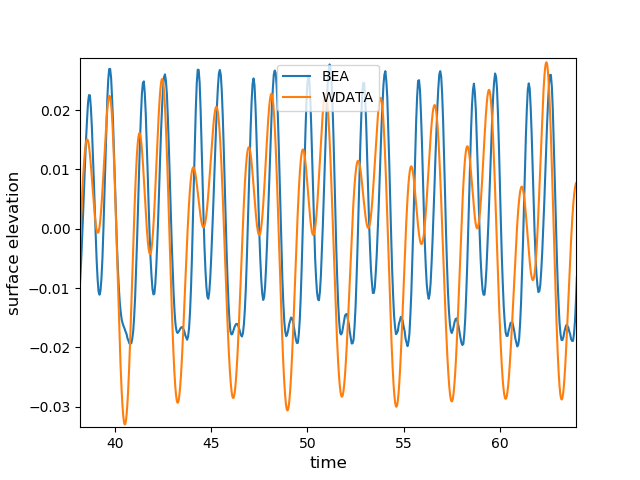}}
\caption{Surface elevation vs time at $x^A_1$ (left) and $x^A_6$ (right).}\label{fig:pt16}
\end{figure}

In Fig. \ref{fig:pt16}, the numerical solution and the experimental measurements at $x_1^A$ and $x_6^A$ are plotted for the different grids.
We note that at $x^A_1$ there is a significant (numerical) dispersion error with $h=0.1$ that is largely gone already with $h=0.05$. This is not surprising since a low-order scheme such as the one used here
requires a rather fine grid to resolve waves accurately. At higher resolutions,
the match with experimental data is very good at $x^A_1$.
\begin{remark}
The experiments were run with a much longer wave train compared to our simulations. Hence, we have cut the plots to the time interval where they overlap.
\end{remark}

At $x^A_6$ we can also observe grid convergence. The difference between the finest and the coarsest is a measure of the numerical errors which appear to be relatively larger for the coarsest mesh than at they were at $x^A_1$.  This is not surprising since the sill induces high-frequency waves that require still higher resolution in space. Furthermore, residual error between the solution at the finest grid at $x^A_6$ and the experimental data is approximately the modelling error. This error is caused by the dispersive relation for the model that does not exactly match the full Euler equations. This is unavoidable, since Boussinesq-type models do not capture all features of the Euler equations. We remark that the dispersive parameters (set 3) have been chosen to produce accurate results for long waves (verified at $x^A_1$) but one could equally well tune them for some other range of frequencies, if that is desirable. 

\begin{remark}
We have verified that the scheme runs stably on a grid with $h=0.001$ but it takes too long to run it all the way to $T=70$ which is why we have omitted it.
\end{remark}

\begin{figure}[h]
\centering
\subfigure{\includegraphics[width=6cm]{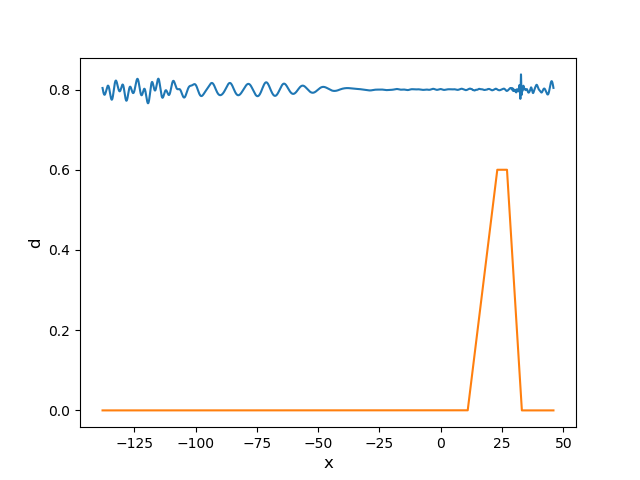}}
\subfigure{\includegraphics[width=6cm]{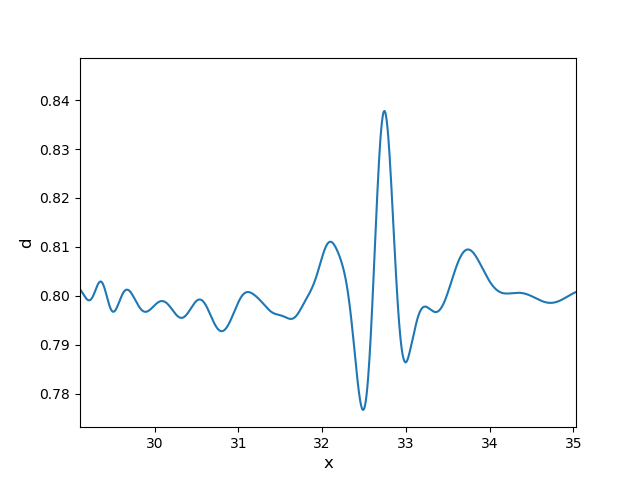}}
\caption{Surface elevation at $T=70$; $h=0.005$.}\label{fig:T70}
\end{figure}

\begin{figure}[h]
\centering
\subfigure[$x^A_2$]{\includegraphics[width=6cm]{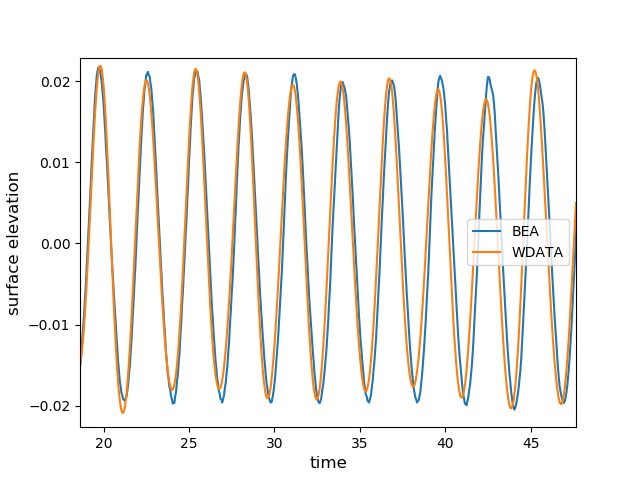}}
\subfigure[$x^A_3$]{\includegraphics[width=6cm]{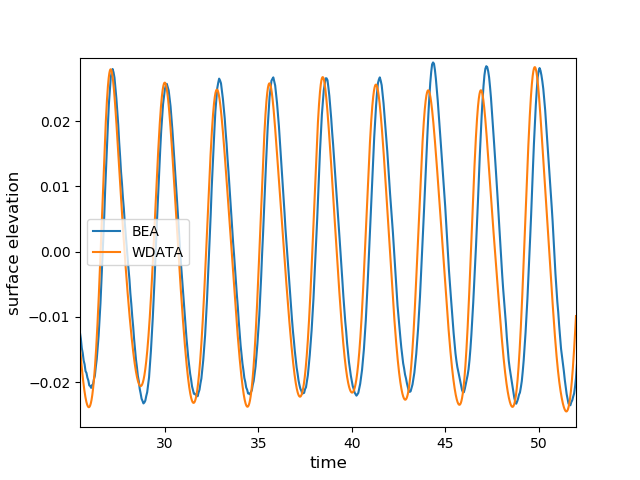}}
\subfigure[$x^A_4$]{\includegraphics[width=6cm]{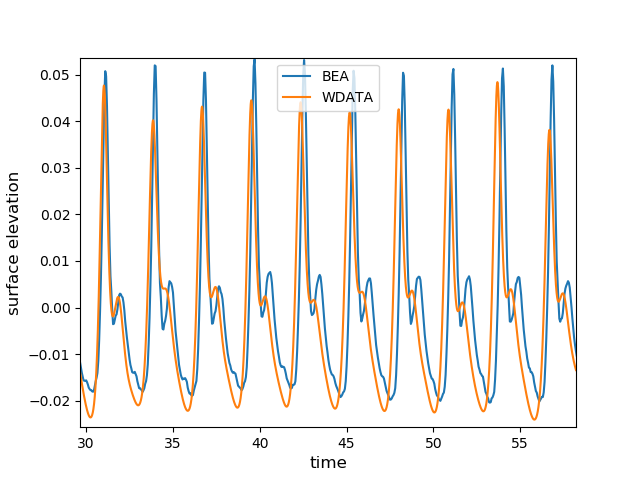}}
\subfigure[$x^A_5$]{\includegraphics[width=6cm]{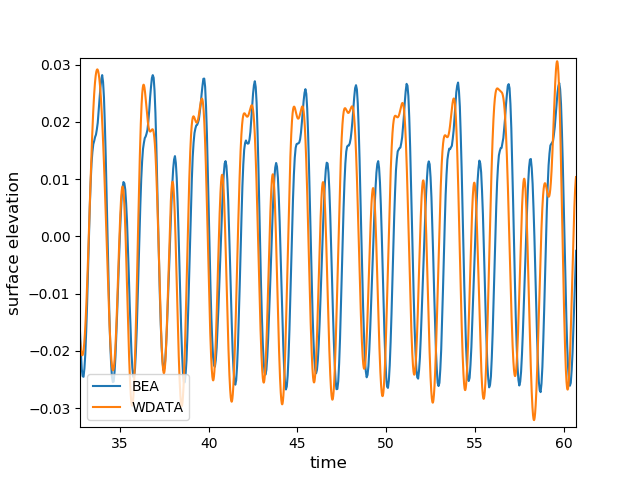}}
\caption{Comparison with experimental data. ($h=0.005$)}\label{fig:xa2345}
\end{figure}

To shed further light on the model, we present the solution for $d$ on the $h=0.005$ grid in Fig. \ref{fig:T70}. We stress that the high-frequency oscillations around $x=32$ are not numerical artefacts but the actual solution. The zoom reveals that the oscillations are fairly well-resolved.
Furthermore, we present a comparison with experimental data at the points $x^A_{2,3,4,5}$ in Fig. \ref{fig:xa2345}. The match with experimental data is fairly good for $x^A_{2,3}$ and less so further downstream. 

\begin{figure}[ht]
  \centering
  \subfigure[$x^A_1$]{\includegraphics[width=6cm]{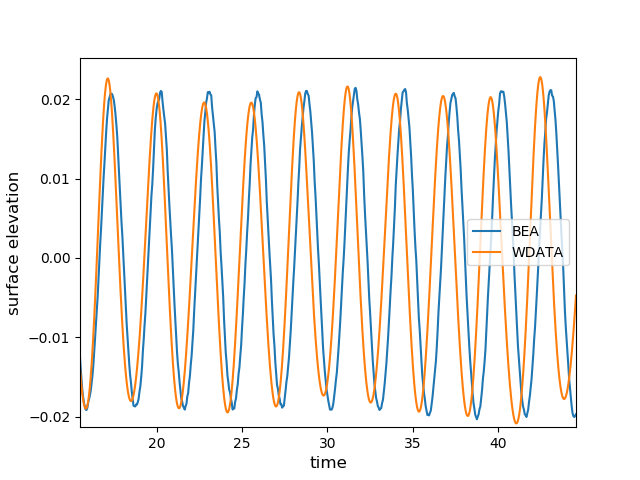}}
\subfigure[$x^A_6$]{\includegraphics[width=6cm]{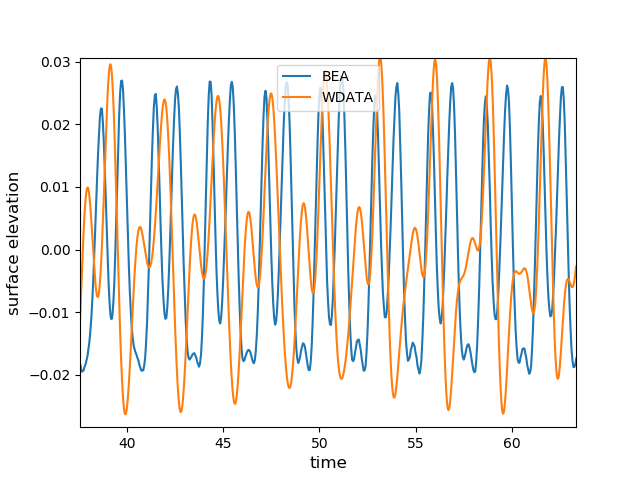}}
\caption{Comparison with experimental data and the parameters ``Set 4'' (\ref{set4})). ($h=0.005$)}\label{fig:set4}
\end{figure}

So far, we have only shown numerical results for ``Set 3'' (\ref{set3}). As discussed above, the model is flexible and can be tuned for particular problems. We end this section by showing in Fig. \ref{fig:set4} the results obtained with $h=0.005$ and ``Set 4'' (\ref{set4}). (We have translated time using the same constant as all previous examples to make comparisons one-to-one.)
For set 4, the dispersion error is larger for long waves than for set 3, which makes the solution less accurate than set 3 at $x^A_1$. The dispersive error is carried downstream making the solution out of phase at $x^A_6$. However, the shape is very similar to the corresponding figure for set 3. (Bottom-right panel in Fig. \ref{fig:pt16}.)

For this problem, ``set 4'' does not produce more accurate solutions than ``set 3''. However, for a problem with only high frequency waves, the situation might be dfferent.

\subsection{Spike}

To demonstrate the robustness of the model and its discrete approximation scheme, we consider a bottom topography with a spike. The domain, the initial data, the CFL number (0.2) and the artificial diffusion coefficient ($\lambda=0.1$) are the same as in section \ref{sec:dinge}. Furthermore, we have arbitrarily used ``Set 4'' (\ref{set4}), since that is of no importance with respect to robustness.

\begin{figure}[h]
\centering
\subfigure[$h=0.01, \,t=0$]{\includegraphics[width=6cm]{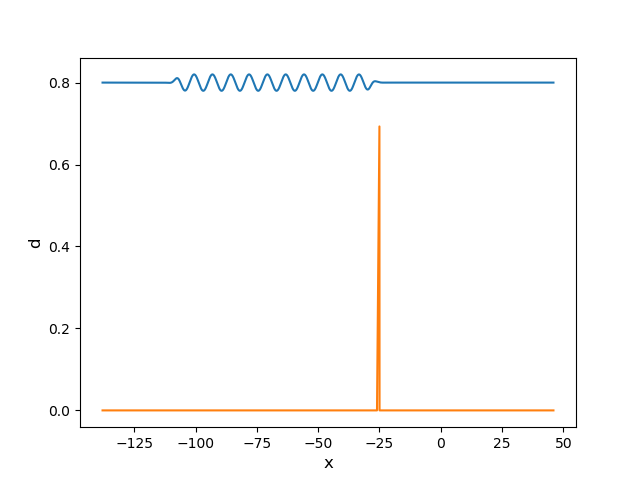}}
\subfigure[$h=0.01, \,t=10$]{\includegraphics[width=6cm]{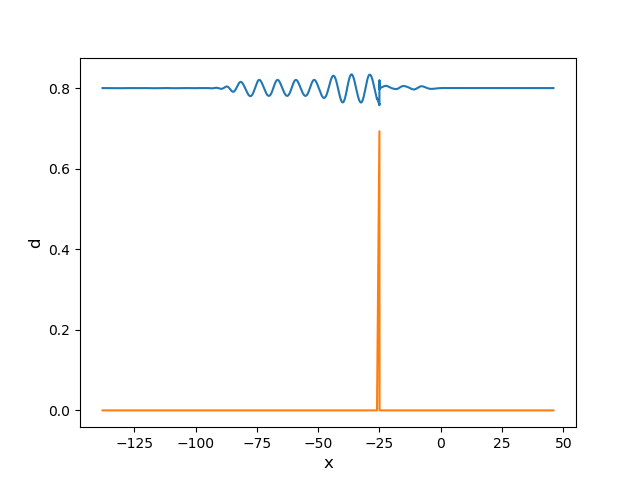}}
\subfigure[$h=0.01,\,t=10$, zoom]{\includegraphics[width=6cm]{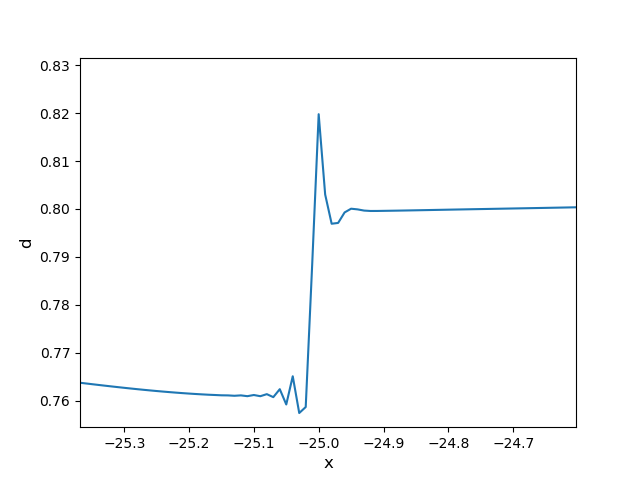}}
\subfigure[$h=0.005,\,t=10$]{\includegraphics[width=6cm]{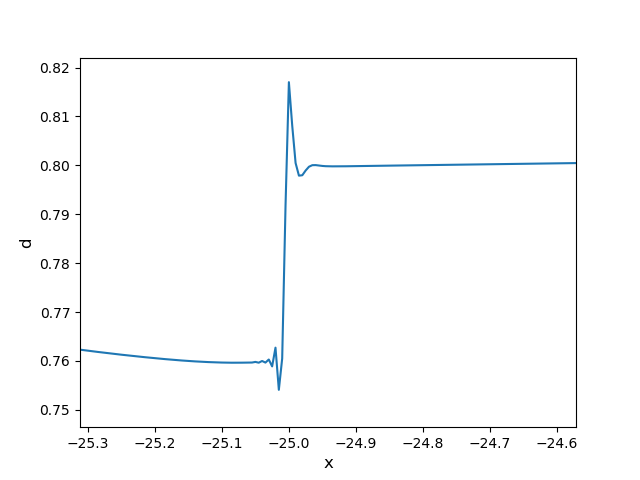}}
\caption{The spike bathymetry}\label{fig:spike}
\end{figure}

The still water depth is 0.8 (as before) apart from at $x=-26...-25$ where a spike is located. It increases linearly to the height $0.7$ in $x\in (-26,-25)$ whereafter there is a jump (true discontinuity) back to zero.  The initial setup for the surface elevation is depicted in Fig. \ref{fig:spike}. (This setup resembles one of the test cases in \cite{LovholtPedersen09}, but the discontinuity makes it even more demanding.)

We run the scheme till $T=5.0$ and do not encounter any stability problems. See Fig. \ref{fig:spike}

\subsection{Cavity}

As another demonstration of robustness, we compute waves travelling over a cavity. We use the same  domain,  initial data, CFL number (0.2) and artificial diffusion coefficient ($\lambda=0.1$) as as in section \ref{sec:dinge}. Also this case was run with ``Set 4'' (\ref{set4}).

\begin{figure}[h]
\centering
\subfigure[$h=0.01,\,t=0$]{\includegraphics[width=6cm]{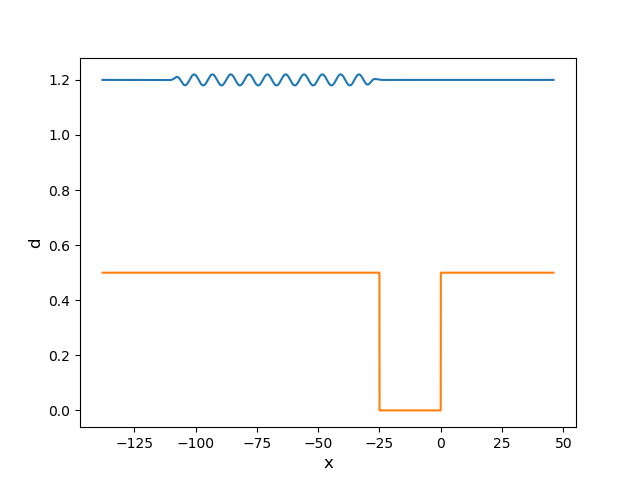}}
\subfigure[$h=0.01, \,t=20$]{\includegraphics[width=6cm]{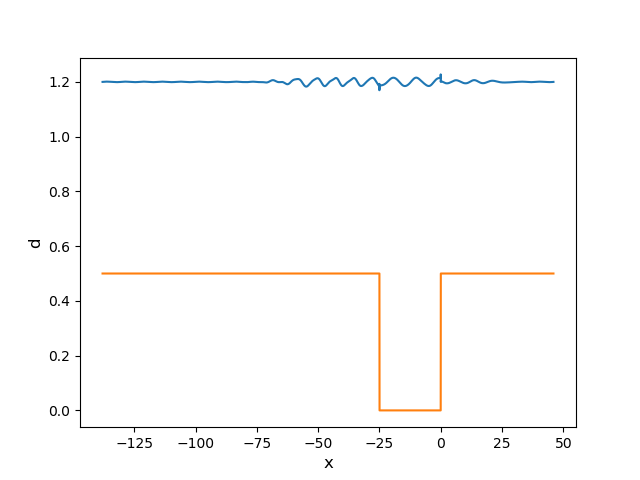}}
\subfigure[$h=0.01,\,t=20$, zoom]{\includegraphics[width=6cm]{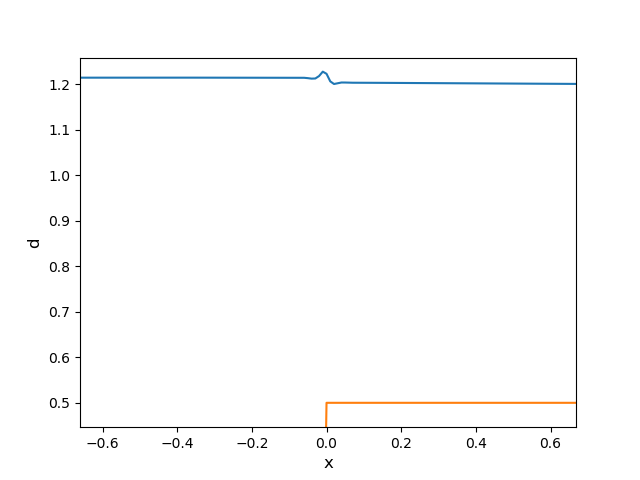}}
\caption{Bathymetry with a cavity}\label{fig:trench}
\end{figure}

The cavity is given by the following bottom topography:
\begin{align*}
  b(x)&=0.5,\quad x<-75\\
  b(x)&=0.0,\quad -75<x<-50\\
  b(x)&=0.5,\quad x>-50\\
\end{align*}
We run the scheme with $h=0.01$. The initial data (for $d$), the bathymetry and the solution at $T=5$, are depicted in Fig. \ref{fig:trench}. 
Clearly, the scheme is stable for this case as well, despite $b(x)$ being discontinuous. Note also the steepening of the waves where the depth is shallower. This is a consequence of the dispersive coefficients being smaller due to the smaller depth leading to a solution that is closer to that of the shallow-water equations.

\section{Two-dimensional extension}

In two space dimensions (2-D), a model should retain the same properties as the one dimensional counterpart, in any arbitrary direction. That is, entropy boundedness and the dispersion relation. 

To this end, we introduce the spatial domain is $(x,y)\in\Omega$, the bathymetry  $b(x,y)$ and the still water depth $h(x,y)$. The dimensionless parameters $\tilde \alpha, \tilde \beta, \tilde \gamma$ take the same values as given by the 1-D analysis (sets 1-4) such that,
\begin{align*}
  \hat \alpha^2&=\tilde\alpha\sqrt{gh(x,y)}h(x,y)^2\\
  \hat \beta &=\tilde \beta h(x,y)^2\\
  \hat \gamma &=\tilde \gamma\sqrt{gh(x,y)}h(x,y)^3
\end{align*}
in complete analogy with the 1-D coefficients. Furthermore, we introuduce the velocity components $v(x,y,t)$ and $w(x,y,t)$ and momentum variables $P=dv$ and $Q=dw$, in the x- and y-directions, respectively.

The 2-D shallow-water equations are rotationally symmetric, and we only need to ensure that the dispersive part is.  We note that for a constant bathymetry, for which the $\tilde \alpha, \tilde \beta,\tilde \gamma$ were obtained to match the dispersion relation of the full model,  the coefficients $\hat\alpha,\hat \beta,\hat \gamma$ are constant. In this case, the dispersive terms associated with $\beta$ and $\gamma$ are linear. Thus, and if $\hat \alpha=0$, the 2-D generalisation of (\ref{eqShallow1})-(\ref{eqShallow2}) with $\tilde \alpha=0$ becomes
\begin{align}
  d_t+P_x+Q_y&= 0\nonumber\\
  P_t + \left(\frac{P^2}{d}\right)_x +(\frac{PQ}{d})_y+ gd(d+b-H)_x &=\nonumber\\
   \left(\hat \beta \left( \hat \beta\frac{P}{d}\right)_{x}\right)_{xt}+\frac{1}{2}\left(\hat \gamma\left(\hat \gamma \frac{P}{d}\right)_x\right)_{xx}+\frac{1}{2}\left(\hat \gamma \left( \frac{P}{d}\right)_{xx}\right)_{x}\label{2dsys}\\
   Q_t + \left(\frac{Q^2}{d}\right)_y +(\frac{PQ}{d})_x+ gd(d+b-H)_y &=\nonumber\\
   \left(\hat \beta \left( \hat \beta\frac{P}{d}\right)_{y}\right)_{yt}+\frac{1}{2}\left(\hat \gamma\left(\hat \gamma \frac{P}{d}\right)_y\right)_{yy}+\frac{1}{2}\left(\hat \gamma \left( \frac{P}{d}\right)_{yy}\right)_{y}\nonumber
\end{align}

\begin{remark}
  Since we have focused on systems where $\hat\alpha=0$, we postpone the generisation of the full system to a future paper. The $\hat\alpha$-terms are non-linear and might not generalise to 2-D in the same straightforward way.
\end{remark}

The system (\ref{2dsys}) is equipped with the entropy (mechanical energy),
\begin{align}
U(\uu)=\frac{1}{2}\left(\frac{P^2}{d} + \frac{Q^2}{d}\right)+\frac{1}{2}gd^2+gdb,\label{U2d}
\end{align}
and the entropy variables, $\ww=\left(g(d+b)-\frac{P^2+Q^2}{2d^2},v,w\right)$. It is straightforward to verify that contracting (\ref{2dsys}) with the entropy variables yields a bound on $U$ in the same way as in the 1-D case.

For brevity, we do not include the 2-D version of the semi-discrete scheme which is straightforward to obtain: The entropy conservative extension of the numerical scheme for the Shallow Water part is found in \cite{Fjordholm_etal11}; the extension of the artificial diffusion terms and the Boussinesq terms is trivial from the 1-D counterparts. The resulting scheme is both etropy stable and well-balanced.

The time discretisation can be done in the same way as in 1-D but will require the solution of a much larger linear system in each time step. This calls for a more sophisticated code and we postpone that work to a future article. 

\section{Conclusions}

Prompted by well-known stability issues with many Boussinesq models, we have derived a new, 
versatile and stable model by taking a novel approach of turning the derivation process around: 
Instead of trying to approximate the Euler system directly, we have obtained
dispersive perturbations to the shallow water system (in 1-D) by requiring both nonlinear stability
and accuracy of the linear dispersion relation.
First, we have choosen the coefficients of the dispersive terms to match the dispersion relation 
of the full Euler equations for a flat bathymetry as accurately as possible with third-order terms.
We have exemplified this technique by proposing a few sets of parameters
that result in systems of different mathematical complexity and different dispersive accuracy. 
These example demonstrate that the system may be tuned to fit specific applications.

Furthermore, we have generalised the system to allow variable bathymetries while retaining 
its entropy-bounded properties. The resulting system reproduces the same dispersive relation 
at any constant depth, and it automatically reverts to the shallow water system 
when the water depth approaches zero. This property dispenses with the necessity 
to turn off the dispersive terms at some small, but otherwise arbitrary, depth; 
a procedure that in many cases introduces instabilities in the computations.

Furthermore, we have derived a nonlinearly entropy-stable semi-discrete scheme for the 
new model and proposed a semi-implicit time-marching scheme. We have demonstrated that the 
scheme is robust for sharp bathymetries, and that it reproduces 
measurements from the Dingemans experiment with reasonable accuracy. 

Finally, a particular version of the system, with dispersive terms added only in the momentum equations 
has been generalised to two spatial dimensions. We have also indicated how the numerical  scheme 
can be generalised to the 2-D system but we postpone a further study in 2-D since 
it requires more sophisticated programming. Future work will also include the analysis 
of entropy-stable boundary conditions and positivity preservation at small depths.

\bibliographystyle{alpha}
\bibliography{ref}



\end{document}